\documentclass[a4paper,11pt]{amsart}
\usepackage[latin1]{inputenc}
\usepackage{amsmath}
\usepackage{amsfonts}
\usepackage{amssymb}
\usepackage{amsthm}
\usepackage{mathrsfs}

\newcommand{\R}{{\mathbb{R}}}

\newcommand{\Z}{{\mathbb{Z}}}
\newcommand{\C}{\mathbb{C}}

\newcommand{\M}{\mathcal{M}}

\def\vecx{{\text{\boldmath$x$}}}
\def\vecy{{\text{\boldmath$y$}}}
\def\vecz{{\text{\boldmath$z$}}}

\def\vecv{{\text{\boldmath$v$}}}
\def\vecm{{\text{\boldmath$m$}}}

\def\vecc{{\text{\boldmath$c$}}}

\def\vecgam{{\text{\boldmath$\gamma$}}}

\def\vecnu{{\text{\boldmath$\nu$}}}
\def\vec0{{\text{\boldmath$0$}}}

\def\Re{\operatorname{Re}}

\newcommand{\ve}{\varepsilon}

\newcommand{\sfrac}[2]{{\textstyle \frac {#1}{#2}}}

\newcommand{\SL}{\mathrm{SL}}

\newtheorem{prop}{Proposition}[section]
\newtheorem{lem}{Lemma}[section]
\newtheorem{thm}{Theorem}

\newtheorem{cor}{Corollary}[section]

\theoremstyle{remark}
\newtheorem{remark}{Remark}

\newenvironment{Proof}{{\bf Proof. }}{\hfill$\square$\newline}
\newenvironment{Proof2}{{\bf Proof of Theorem 2. }}{\hfill$\square$\newline}
\newenvironment{Proof1}{{\bf Proof of Theorem 1. }}{\hfill$\square$\newline}

\begin{document}
\title[On the value distribution and moments of the Epstein zeta function]{On the value distribution and moments of the Epstein zeta function to the right of the critical strip}
\author{Anders S\"odergren}
\address{Department of Mathematics, Uppsala University, Box 480,\newline
\rule[0ex]{0ex}{0ex}\hspace{8pt} SE-75106 Uppsala, Sweden\newline
\rule[0ex]{0ex}{0ex}\hspace{8pt} {\tt sodergren@math.uu.se}}
\date{6 June 2010}

\maketitle

\begin{abstract}
We study the Epstein zeta function $E_n(L,s)$ for $s>\frac{n}{2}$ and determine for fixed $c>\frac{1}{2}$ the value distribution and moments of $E_n(\cdot,cn)$ (suitably normalized) as $n\to\infty$. We further discuss the random function $c\mapsto E_n(\cdot,cn)$ for $c\in[A,B]$ with $\frac{1}{2}<A<B$ and determine its limit distribution as $n\to\infty$.
\end{abstract}

\tableofcontents


\section{Introduction}

Let $G_n=\SL(n,\R)$ and $\Gamma_n=\SL(n,\Z)$. We will be interested in the space $X_n=\Gamma_n\backslash G_n$ considered as the space of $n$-dimensional lattices of covolume $1$. Here $\Gamma_ng$ corresponds to the lattice $\Z^ng\subset\R^n$. We let $\mu_n$ denote the Haar measure on $G_n$, normalized so that it represents the unique right $G_n$-invariant probability measure on the homogeneous space $X_n$.

For $L\in X_n$ and $\Re s>\frac{n}{2}$ the Epstein zeta function is defined by
\begin{align}\label{Epstein}
E_n(L,s)={\sum_{\vecm\in L}}'|\vecm|^{-2s},
\end{align}
where $'$ denotes that the zero vector should be omitted. $E_n(L,s)$ has an analytic continuation to $\C$ except for a simple pole at $s=\frac{n}{2}$ and furthermore it satisfies the functional equation
\begin{align*}
F_n(L,s)=F_n(L^*,\sfrac{n}{2}-s),
\end{align*}
where $F_n(L,s):=\pi^{-s}\Gamma(s)E_n(L,s)$ and $L^*$ is the dual lattice of $L$.

The Epstein zeta function is in many ways analogous to the Riemann zeta function. In particular we have the relation
\begin{align*}
 \zeta(2s)=\frac{1}{2}E_1(\Z,s).
\end{align*}
With this analogy in mind it is natural to call the region $0<\Re s<\frac{n}{2}$ the critical strip for $E_n(L,s)$. Note however that for all $n\geq2$ there exist lattices $L\in X_n$ for which the Riemann hypothesis for $E_n(L,s)$ is known to fail (cf. \cite[Thm.\ 1]{terras3}; see also \cite{bateman}, \cite{stark}, \cite{terras2} and \cite{terras}).

In this paper we will study $E_n(L,s)$ with $n$ large and $s>\frac{n}{2}$. In particular we will, for fixed $c>\frac{1}{2}$, be interested in questions concerning the value distribution of $E_n(L,cn)$ as $n\to\infty$. These questions are mainly motivated by the work of Sarnak and Strömbergsson (\cite[Sec.\ 6]{sst}) on the height function for flat tori in large dimensions. For the flat torus $\R^n/L$, with $L\in X_n$, the height function is defined by
\begin{align*}
 h_n(\R^n/L)=2\log 2\pi+\frac{\partial}{\partial s}E_n(L^*,s)_{|s=0}.
\end{align*}
Sarnak and Strömbergsson show that the function $h_n$ concentrates at the value $\log 4\pi-\gamma+1$ \footnote{Here $\gamma$ is Euler's constant.} as $n\to\infty$. Interpreted in terms of the Epstein zeta function this result states that if $\ve>0$ is fixed then
\begin{align*}
\text{Prob}_{\mu_n}\Big\{L\in X_n\,\,\Big|\,\,\big|\sfrac{\partial}{\partial s}E_n(L,s)_{|s=0}-(1-\gamma-\log \pi)\big|<\ve\Big\}\to 1
\end{align*}
as $n\to\infty$. (In this connection, recall that $E_n(L,0)=-1$ for all $n$ and all $L\in X_n$.)

This description of the derivative of $E_n(L,s)$ at the point $s=0$ naturally suggests the question of what we can say about the value distribution of the function $E_n(L,s)$ itself as $n\to\infty$. We start this investigation to the right of the critical strip since there $E_n(L,s)$ is given by the simple formula \eqref{Epstein}. The study of $E_n(L,s)$ for large $s$ is further related to the classical problem of finding the densest lattice sphere packing in $\R^n$. To be more precise, Ryshkov (\cite{ryshkov}) has showed that the densest lattice sphere packing of $\R^n$ is given by the lattice that minimizes $E_n(L,s)$ as $s\to\infty$.

In the present context it is most natural to consider a normalized version of $E_n(L,s)$. Our main result is that the value distribution, as $n\to\infty$, of the normalized Epstein zeta function can be completely described in terms of the points of a Poisson process on the positive real line.

\begin{thm}\label{mainthm}
Let $V_n$ denote the volume of the $n$-dimensional unit ball. Let $\mathcal P$ be a Poisson process on the positive real line with intensity $\frac{1}{2}$ and let $T_1,T_2,T_3,\ldots$ denote the points of $\mathcal P$ ordered so that $0<T_1< T_2< T_3<\cdots$. Then, for fixed $c>\frac{1}{2}$, the distribution of the random variable $V_n^{-2c}E_n(\cdot,cn)$ converges to the distribution of $2\sum_{j=1}^{\infty}T_j^{-2c}$ as $n\to\infty$. In fact, for any $m\geq1$ and fixed $\frac{1}{2}<c_1<\cdots<c_m$, the distribution of the random vector
\begin{align*}
\Big(V_n^{-2c_1}E_n(\cdot,c_1n),\ldots,V_n^{-2c_m}E_n(\cdot,c_mn)\Big)
\end{align*}
converges to the distribution of
\begin{align*}
\Big(2\sum_{j=1}^{\infty}T_j^{-2c_1},\ldots,2\sum_{j=1}^{\infty}T_j^{-2c_m}\Big)
\end{align*}
as $n\to\infty$.
\end{thm}

Actually, with the same notation as in Theorem \ref{mainthm}, even more is true:

\begin{thm}\label{curvethm}
For each $n\in\Z^+$ and any fixed $\frac{1}{2}<A<B$ consider
\begin{align*}
 c\mapsto V_n^{-2c}E_n(\cdot,cn)
\end{align*}
as a random function in $C[A,B]$. The distribution of this random function converges to the distribution of
\begin{align*}
 c\mapsto 2\sum_{j=1}^{\infty}T_j^{-2c}
\end{align*}
as $n\to\infty$.
\end{thm}

The most important ingredient in the proof of Theorem \ref{mainthm} is our result \cite{jag} on the distribution of lengths of lattice vectors in a random lattice in $X_n$. It states that, as $n\to\infty$, the suitably normalized non-zero vector lengths in a random lattice $L\in X_n$ behave like the points of a Poisson process on the positive real line. This result together with equation \eqref{Epstein} suggest that the limit distribution of the normalized Epstein zeta function should be that of the "random Dirichlet series" appearing in Theorem \ref{mainthm} and Theorem \ref{curvethm}. It is immediate from \cite{jag} that the ``symmetric'' partial sums of $V_n^{-2c}E_n(\cdot,cn)$ converge in distribution to the corresponding partial sums of $2\sum_{j=1}^{\infty}T_j^{-2c}$ as $n\to\infty$. The proof is finished by an approximation argument using a bound on the second moment.

We can strengthen the result in Theorem \ref{mainthm} by showing that we do not only have convergence in distribution but also convergence in moments after an explicit and tractable truncation. (A truncation is necessary already in order for the moments of $E_n(\cdot,s)$ to exist.) We will consider the truncation $E_{R_n(\delta)}(L,s)$ of $E_n(L,s)$ that discards the contribution to $E_n(L,s)$ from all lattice vectors in $L$ belonging to the $n$-ball of volume $\delta>0$ centered at the origin:
\begin{align*}
E_{R_n(\delta)}(L,s)=\sum_{\substack{\vecm\in L\\|\vecm|>R_n(\delta)}}|\vecm|^{-2s},\qquad
R_n(\delta)=\Big(\frac{\delta}{V_n}\Big)^{\frac{1}{n}}. 
\end{align*}
It follows that $E_n(L,s)-E_{R_n(\delta)}(L,s)$ is non-zero only on a set of measure at most $\frac{1}{2}\delta$ (cf.\ \eqref{secondalt} below). We show that the moments of $V_n^{-2c}E_{R_n(\delta)}(\cdot,cn)$ converge to those of a similar truncation of $2\sum_{j=1}^{\infty}T_j^{-2c}$ as $n\to\infty$. In precise terms:

\begin{thm}\label{momentthm}
Let $c>\frac{1}{2}$  and $\delta>0$ be fixed. Let further $I(\cdot)$ be the indicator function. Then every moment of $V_n^{-2c}E_{R_n(\delta)}(\cdot,cn)$ converges to the corresponding moment of $2\sum_{j=1}^{\infty}I\big(T_j>\delta\big)T_j^{-2c}$ as $n\to\infty$. Furthermore, for any $m\geq1$ and fixed $\frac{1}{2}<c_1<\cdots<c_m$, the joint moments of the random vector
\begin{align*}
\Big(V_n^{-2c_1}E_{R_n(\delta)}(\cdot,c_1n),\ldots,V_n^{-2c_m}E_{R_n(\delta)}(\cdot,c_mn)\Big)
\end{align*}
converge to the corresponding joint moments of
\begin{align*}
\Big(2\sum_{j=1}^{\infty}I\big(T_j>\delta\big)T_j^{-2c_1},\ldots,2\sum_{j=1}^{\infty}I\big(T_j>\delta\big)T_j^{-2c_m}\Big)
\end{align*}
as $n\to\infty$.
\end{thm}

The proof of Theorem \ref{momentthm} is based on the approach by Sarnak and Strömbergsson. We determine the limiting expressions for the moments of $V_n^{-2c}E_{R_n(\delta)}(\cdot,cn)$ using the integration formula of Siegel (\cite{siegel}) and its generalization by Rogers (\cite{rogers1}). The limits of the mean and variance of $V_n^{-2c}E_{R_n(\delta)}(\cdot,cn)$ are found by straightforward calculations. When $k\geq3$ the $k$:th moment requires more advanced estimates in order to find the limiting value. These estimates are based on a mixture of methods discussed by Rogers in \cite{rogers2}, \cite{rogers3} and \cite{rogers4}. The exact formula for the limit of the $k$:th order moment is given in
Theorem \ref{higher} in Section \ref{moments}.

Theorem \ref{momentthm} implies Theorem \ref{mainthm} by letting $\delta\to0$, cf.\ Section \ref{altsec}. We expect that this approach to the value distribution of $E_n(L,s)$ can be used also to study the Epstein zeta function in the critical strip. The author plans to address this more difficult problem in a forthcoming paper.

\section{Preliminaries}

\subsection{Siegel's and Rogers' integration formulae}

In this section we fix some notation concerning the integration formula that will be the major technical tool when discussing the moments of $E_{R_n(\delta)}(L,s)$.

Let $f:\R^n\to\R$ be a non-negative Borel measurable function. In \cite{siegel} Siegel proves the mean value formula
\begin{align}\label{siegelform}
\int_{X_n}\sum_{\vecm\in L\setminus\{\vec0\}}f(\vecm)\,d\mu_n(L)=\int_{\R^n}f(\vecx)\,d\vecx.
\end{align}
We next describe Rogers' generalization of Siegel's formula.

Let $1\leq k\leq n-1$ and let $\rho:(\R^n)^k\to\R$ be a non-negative Borel measurable function. In \cite{rogers1} Rogers considers the integral
\begin{align*}
\int_{X_n}\sum_{\vecm_1,\ldots,\vecm_k\in L}\rho(\vecm_1,\ldots,\vecm_k)\,d\mu_n(L),
\end{align*}
and shows that it equals a certain (positive) infinite linear combination of integrals of $\rho$ over various linear subspaces of $(\R^n)^k$. In this paper we will be interested in the similar integral
\begin{align}\label{rogint}
\int_{X_n}\sum_{\vecm_1,\ldots,\vecm_k\in L\setminus\{\vec0\}}\rho(\vecm_1,\ldots,\vecm_k)\,d\mu_n(L).
\end{align}
It follows from Rogers' formula in \cite{rogers1} that the integral in \eqref{rogint} equals
\begin{align}\label{rogform}
&\int_{\R^n}\cdots\int_{\R^n}\rho(\vecx_1,\ldots,\vecx_k)\,d\vecx_1\ldots d\vecx_k\nonumber\\
&+\sum_{(\nu,\mu)}\sum_{q=1}^{\infty}\sum_{D}\Big(\frac{e_1}{q}\cdots\frac{e_m}{q}\Big)^n\int_{\R^n}\cdots\int_{\R^n}\rho\Big(\sum_{i=1}^m\frac{d_{i1}}{q}\vecx_i,\ldots,\sum_{i=1}^m\frac{d_{ik}}{q}\vecx_i\Big)\,d\vecx_1\ldots d\vecx_m.
\end{align}
Here the outer sum is over all divisions $(\nu,\mu)=(\nu_1,\ldots,\nu_m;\mu_1,\ldots,\mu_{k-m})$ of the numbers $1,\ldots,k$ into two sequences $\nu_1,\ldots,\nu_m$ and $\mu_1,\ldots,\mu_{k-m}$ with $1\leq m\leq k-1$, satisfying
\begin{align}\label{division}
& 1\leq \nu_1<\nu_2<\ldots<\nu_m\leq k,\nonumber\\
& 1\leq\mu_1<\mu_2<\ldots<\mu_{k-m}\leq k,\\
& \nu_i\neq\mu_j, \text{ if $1\leq i\leq m$, $1\leq j \leq k-m$}.\nonumber
\end{align}
The inner sum in \eqref{rogform} is over all $m\times k$ matrices $D$, with no column vanishing, with integer elements having greatest common divisor equal to $1$, and with
\begin{align*}
&d_{i\nu_j}=q\delta_{ij}, \hspace{10pt}i=1,\ldots,m,\,\,j=1,\ldots,m,\nonumber\\
&d_{i\mu_j}=0,\hspace{10pt}\text{ if }\,\mu_j<\nu_i,\,\,i=1,\ldots,m,\,\,j=1,\ldots,k-m.
\end{align*}
We call these matrices $(\nu,\mu)$-admissible. A matrix is called $k$-admissible if it is $(\nu,\mu)$-admissible for some division $(\nu,\mu)$ satisfying \eqref{division}. Finally $e_i=(\ve_i,q)$, $i=1,\ldots,m$, where $\ve_1,\ldots,\ve_m$ are the elementary divisors of the matrix $D$.


\begin{remark}
It follows from the conditions on the matrices $D$ above and \cite[Thm.\ 14.5.1]{hua} that in all cases we have $e_1=1$. In particular it follows that we always have
\begin{align}\label{element}
\Big(\frac{e_1}{q}\cdots\frac{e_m}{q}\Big)^n\leq q^{-n}.
\end{align}
\end{remark}

\subsection{Recapitulation of results in \cite{jag}}\label{Latticepoisson}

Given a lattice $L\in X_n$, we order the non-zero vector lengths in $L$ as $0<\ell_1\leq \ell_2\leq \ell_3\leq\ldots$, where we count the common length of the vectors $\vecx$ and $-\vecx$ only once. For $j\geq1$, we define
\begin{align}\label{lengthvolume}
 \mathcal V_j(L):=\frac{\pi^{n/2}}{\Gamma(\frac{n}{2}+1)}\ell_j^n
\end{align}
so that $\mathcal V_j(L)$ is the volume of an $n$-dimensional ball of radius $\ell_j$. Finally, for $t\geq0$, we let
\begin{align*}
 \tilde N_{t}(L):=\#\{j:\mathcal{V}_j(L)\leq t\}.
\end{align*}
The main theorem in \cite{jag} is the following:

\begin{thm}\label{thmjag}
Let $\{N(t),t\geq0\}$ be a Poisson process on the positive real line with intensity $\frac{1}{2}$. Then the stochastic process $\{\tilde N_{t}(\cdot),t\geq0\}$ converges weakly to \mbox{$\{N(t),t\geq0\}$} as $n\to\infty$.
\end{thm}

Let us consider the Poisson process $\{N(t),t\geq0\}$. We recall that $N(t)$ denotes the number of points falling in the interval $(0,t]$ and that $N(t)$ is Poisson distributed with expectation value $\frac{1}{2}t$. We let $T_1,T_2,T_3,\ldots$ denote the points of the process ordered in such a way that $0<T_1<T_2< T_3<\cdots$. 

The proof of Theorem \ref{thmjag} is based on calculations using Rogers' integration formula and the following two observations, which we recall here for convenience.

\begin{prop}\label{svlem}
 Let $k\geq1$ and denote by $\mathcal{P}(k)$ the set of partitions of $\{1,\ldots,k\}$.  For $1\leq j\leq k$ let $f_j:\R_{\geq0}\to\R$ be functions satisfying $\prod_{j\in B}f_j\in L^1(\R_{\geq0})$ for every nonempty subset $B\subseteq\{1,\ldots,k\}$. Then
\begin{align*}
\mathbb E\Big(\prod_{j=1}^k\Big(\sum_{n=1}^{\infty}f_j(T_n)\Big)\Big)=\underset{P \in\mathcal{P}(k)}{\sum}2^{-\#P}\underset{B\in P}{\prod}\Big(\int_{0}^{\infty}\prod_{j\in B}f_j(x)\,dx\Big).
\end{align*}
\end{prop}

\begin{lem}\label{bijection}
Let $k\geq1$ and let $\mathcal{D}(k)$ be the set of $k$-admissible matrices $D$ with all entries in $\{0,1\}$ and exactly one entry equal to $1$ in each column, together with the $k\times k$ indentity matrix $I_k$. When $D=I_k$ let $\nu_i=i$, $1\leq i\leq k$. Then there is a bijection $g:\mathcal{D}(k)\to\mathcal{P}(k)$ with the property that if $D\in\mathcal{D}(k)$ is an $m\times k$ matrix and $g(D)=P=\{B_1,\ldots,B_{\#P}\}$ then $\#P=m$ and $\{\nu_1,\ldots,\nu_m\}=\{\min_{j\in B_1} j,\ldots,\min_{j\in B_m} j\}$.
\end{lem}

\subsection{Normalization of $E_n(L,s)$}

For $n\geq1$ we let
\begin{align*}
 V_n:=\frac{\pi^{n/2}}{\Gamma(\frac{n}{2}+1)}
\end{align*}
so that $V_n$ is the volume of the unit ball in $\R^n$. For future reference we also recall that
\begin{align}\label{unitvolume}
 V_n=\frac{\omega_n}{n},
\end{align}
where $\omega_n$ is the volume of the $n-1$ sphere. For $s>\frac{n}{2}$ we can, using \eqref{lengthvolume}, write $E_n(L,s)$ in the form
\begin{align*}
E_{n}(L,s)={\sum_{\vecm\in L}}'|\vecm|^{-2s}
=V_n^{\frac{2s}{n}}{\sum_{\vecm\in L}}'\big(V_n|\vecm|^n\big)^{-\frac{2s}{n}}
=2V_n^{\frac{2s}{n}}\sum_{j=1}^{\infty}\mathcal{V}_j(L)^{-\frac{2s}{n}}.
\end{align*}
Hence it is natural to consider the normalized function
\begin{align*}
\mathcal{E}_{n}(L,s):=V_n^{\frac{-2s}{n}}E_n(L,s),
\end{align*}
so that 
\begin{align}\label{norm2}
\mathcal{E}_{n}(L,s)=2\sum_{j=1}^{\infty}\mathcal{V}_j(L)^{-\frac{2s}{n}}
\end{align}
for $s>\frac{n}{2}$. It is this normalized form of the Epstein zeta function that will be in focus in the present paper.

\subsection{Truncation of $E_n(L,s)$}

By an application of Rogers' formula \eqref{rogform} it is clear that the moments of $E_{n}(\cdot,s)$ do not exist for any $s>\frac{n}{2}$. In order to discuss moments we will thus need to consider a truncation of the Epstein zeta function. Here we will focus on the truncation $E_{R_n(\delta)}(L,s)$ of $E_n(L,s)$ that discards the contribution to $E_n(L,s)$ from all lattice vectors in $L$ belonging to the $n$-ball of volume $\delta>0$ centered at the origin. The details are as follows:

For $n\geq1$ and $\delta>0$ we let $R_n(\delta)$ to be the radius of a ball of volume $\delta$ in $\R^n$. Hence
\begin{align}\label{rnd}
R_n(\delta)=\Big(\frac{\delta}{V_n}\Big)^{\frac{1}{n}},
\end{align}
and for $\vecx\in\R^n$ we note that
\begin{align}\label{volume}
|\vecx|>R_n(\delta)\Leftrightarrow V_n|\vecx|^n>\delta.
\end{align}
As is indicated above we define
\begin{align}\label{cutoff}
E_{R_n(\delta)}(L,s):={\sum_{\vecm\in L}}|\vecm|^{-2s}I_{R_n(\delta)}(\vecm),
\end{align}
where $I_{R_n(\delta)}$ is the cut-off function
\begin{align*}
 I_{R_n(\delta)}(\vecx)=
\begin{cases}
1 & \text{ if $|\vecx|>R_n(\delta)$}\\
0 & \text{ if $|\vecx|\leq R_n(\delta)$.}
\end{cases}
\end{align*}

Applying Siegel's formula \eqref{siegelform} to the function $\chi_{R_n(\delta)}(\vecx)=1-I_{R_n(\delta)}(\vecx)$ we find that 
\begin{align}\label{secondalt}
\int_{X_n}\sum_{\vecm\in L\setminus\{\vec0\}}\chi_{R_n(\delta)}(\vecm)\,d\mu_n(L)=\int_{\R^n}\chi_{R_n(\delta)}(\vecx)\,d\vecx=\delta.
\end{align}
It follows that $E_n(L,s)-E_{R_n(\delta)}(L,s)$ is non-zero on a set of measure at most $\frac{1}{2}\delta$.

Also when discussing moments of the truncated Epstein zeta function we will find it most natural to work with the normalized form   
\begin{align}\label{norm3}
\mathcal{E}_{R_n(\delta)}(L,s):=V_n^{\frac{-2s}{n}}E_{R_n(\delta)}(L,s)=2\sum_{\mathcal{V}_j(L)>\delta}\mathcal{V}_j(L)^{-\frac{2s}{n}}.
\end{align}

\subsection{The random variable $T(c)$}

Let us return to the Poisson process $\{N(t),t\geq0\}$ on the positive real line with constant intensity $\frac{1}{2}$. As in Section \ref{Latticepoisson} we let $T_j$, $j\geq1$, be the points of the process taken in increasing order. Let $c>\frac{1}{2}$. Motivated by Theorem \ref{thmjag} and equation \eqref{norm2} we will be interested in the random variable
\begin{align*}
T(c):=2\sum_{j=1}^{\infty}T_j^{-2c}.
\end{align*}
The purpose of this paper is to prove that $T(c)$ is (in several senses) the limit of $\mathcal{E}_{n}(\cdot,cn)$ as $n\to\infty$.

Using Proposition \ref{svlem} it is easy to see that the moments of $T(c)$ do not exist. To make it possible to anyhow enter a discussion of moments we introduce the truncated version 
\begin{align*}
T(c,\delta):=2\sum_{j=1}^{\infty}I\big(T_j>\delta\big)T_j^{-2c},
\end{align*}
where $\delta>0$ and $I(\cdot)$ is the indicator function. We remark that this truncation is chosen to match the truncation $\mathcal{E}_{R_n(\delta)}(L,s)$ of $\mathcal{E}_{n}(L,s)$ introduced in the previous section.

We now calculate the moments of $T(c,\delta)$. From Proposition \ref{svlem}, with $k\geq1$ and $f_1(x)=\ldots=f_k(x)=2I(x>\delta)x^{-2c}$, we get
\begin{align}\label{Tmoment}
\mathbb E\big(T(c,\delta)^k\big)&=2^k\mathbb E\Big(\Big(\sum_{j=1}^{\infty}I\big(T_j>\delta\big)T_j^{-2c}\Big)^k\Big)\nonumber\\
&=2^k\underset{P \in\mathcal{P}(k)}{\sum}2^{-\#P}\underset{B\in P}{\prod}\int_{\delta}^{\infty}x^{-2c\#B}\,dx\nonumber\\
&=2^k\underset{P \in\mathcal{P}(k)}{\sum}2^{-\#P}\underset{B\in P}{\prod}\frac{\delta^{1-2c\#B}}{2c\#B-1}\nonumber\\
&=2^k\underset{P \in\mathcal{P}(k)}{\sum}2^{-\#P}\delta^{\#P-2kc}\underset{B\in P}{\prod}\frac{1}{2c\#B-1}.
\end{align}

It is now easy to show that the distribution of $T(c,\delta)$ is uniquely determined by its moments.

\begin{lem}\label{determ}
Let $c>\frac{1}{2}$  and $\delta>0$ be fixed. Let $\mu_{T(c,\delta)}$ be the distribution of $T(c,\delta)$. Then $\mu_{T(c,\delta)}$ is the only probability measure on the positive real line with the moments $\mathbb E\big(T(c,\delta)^k\big)$, $k\geq1$.
\end{lem}

\begin{Proof}
By \cite[Thm.\ 30.1]{billing} it is enough to prove that the moment generating function, $\psi_{T(c,\delta)}(t)$, has a positive radius of convergence. We recall that
\begin{align*}
\psi_{T(c,\delta)}(t)=\sum_{k=0}^{\infty}\frac{\mathbb E\big(T(c,\delta)^k\big)}{k!}t^k.
\end{align*}
We call the coefficients in this power series $\alpha_k$ and study the radius of convergence using the ratio test. By dividing $\mathcal{P}(k+1)$ into the partitions that have $\{k+1\}$ as an element and those who do not have $\{k+1\}$ as an element, we find that
\begin{align}\label{carl}
 \frac{\alpha_{k+1}}{\alpha_k}&=\frac{1}{k+1}\frac{\mathbb E\big(T(c,\delta)^{k+1}\big)}{\mathbb E\big(T(c,\delta)^k\big)}
=\frac{2}{k+1}\frac{\underset{P \in\mathcal{P}(k+1)}{\sum}2^{-\#P}\delta^{\#P-2(k+1)c}\underset{B\in P}{\prod}\frac{1}{2c\#B-1}}{\underset{P' \in\mathcal{P}(k)}{\sum}2^{-\#P'}\delta^{\#P'-2kc}\underset{B'\in P'}{\prod}\frac{1}{2c\#B'-1}}\nonumber\\
&\leq\frac{2}{k+1}\Bigg(\frac{\delta^{1-2c}}{2(2c-1)}+\frac{\underset{P \in\mathcal{P}(k)}{\sum}2^{-\#P}\delta^{\#P-2(k+1)c}\#P\underset{B\in P}{\prod}\frac{1}{2c\#B-1}}{\underset{P' \in\mathcal{P}(k)}{\sum}2^{-\#P'}\delta^{\#P'-2kc}\underset{B'\in P'}{\prod}\frac{1}{2c\#B'-1}}\Bigg)\nonumber\\
&\leq\frac{2}{k+1}\bigg(\frac{\delta^{1-2c}}{2(2c-1)}+k\delta^{-2c}\bigg)\leq M,
\end{align}
for all $k\geq1$ and some $0<M<\infty$ (depending on $c$ and $\delta$). It follows that $\psi_{T(c,\delta)}(t)$ has a positive (or infinite) radius of convergence.
\end{Proof}

We end this section by discussing two random vectors related to $T(c,\delta)$. Let $m\geq1$ and fix $\vecc=(c_1,\ldots,c_m)$ satisfying $\frac{1}{2}<c_1< c_2<\cdots< c_m$. In Theorem \ref{mainthm} and Theorem \ref{momentthm} we will be  interested in
\begin{align}\label{Tvec}
T(\vecc):=\big(T(c_1),\ldots,T(c_m)\big)
\end{align}
and
\begin{align}\label{Tdvec}
T(\vecc,\delta):=\big(T(c_1,\delta),\ldots,T(c_m,\delta)\big)
\end{align}
respectively. As before the truncation $T(\vecc,\delta)$ is introduced in order to discuss moments.
To be more precise we let $\kappa\geq1$ and fix $\gamma_j\in\{c_1,\ldots,c_m\}$, $1\leq j\leq\kappa$, satisfying $\gamma_1\leq \gamma_2\leq\cdots\leq \gamma_{\kappa}$. Using Proposition \ref{svlem}, with $f_j(x)=2I(x>\delta)x^{-2\gamma_j}$ ($1\leq j\leq \kappa$), we get
\begin{align}\label{Tmixedmoment}
\mathbb E\Big(\prod_{j=1}^{\kappa}T(\gamma_j,\delta)\Big)&=\mathbb E\Big(\prod_{j=1}^{\kappa}\Big(\sum_{n=1}^{\infty}f_j(T_n)\Big)\Big)\nonumber\\
&=2^{\kappa}\underset{P \in\mathcal{P}(\kappa)}{\sum}2^{-\#P}\delta^{\#P-2\sum_{j=1}^{\kappa}\gamma_j}\underset{B\in P}{\prod}\frac{1}{2B_{\vecgam}-1},
\end{align}
where, for $B\in P\in\mathcal{P}(\kappa)$, we have $B_{\vecgam}=\sum_{j\in B}\gamma_j$.

\begin{lem}\label{multideterm}
Let $\delta>0$ be fixed. Let $m\geq1$ and fix $\vecc=(c_1,\ldots,c_m)$ satisfying $\frac{1}{2}<c_1< c_2<\cdots< c_m$. Let $\mu_{T(\vecc,\delta)}$ be the distribution of the random vector $T(\vecc,\delta)$. Then $\mu_{T(\vecc,\delta)}$ is the only probability measure on $\R^m$ with joint moments given by \eqref{Tmixedmoment}.
\end{lem}

\begin{Proof}
We use \cite[Thm.\ 2.3]{dejeu} with $\{\vecv_1,\ldots,\vecv_{m}\}$ equal to the standard basis in $\R^{m}$. It follows that it is enough to show that the sequences $\{\mathbb E\big(T(c_j,\delta)^k\big)\}_{k=1}^{\infty}$, $1\leq j\leq m$, satisfy the Carleman condition, i.e.\ that
\begin{align*}
 \sum_{k=1}^{\infty}\frac{1}{\mathbb E\big(T(c_j,\delta)^{2k}\big)^{1/2k}}=\infty,\hspace{30pt}1\leq j\leq m.
\end{align*}
But the divergence of these series follows from \eqref{carl} and an application of Stirling's formula.
\end{Proof}

\section{The moments of the Epstein zeta function}\label{moments}

In this section we will discuss the moments of the truncated and normalized Epstein zeta function $\mathcal{E}_{R_n(\delta)}(\cdot,cn)$ for fixed $c>\frac{1}{2}$. For each fixed moment we will assume that $n$ is large enough to make Rogers' integration formula applicable. The goal is to study the moments as $n\to\infty$.

As an easy first step we calculate the expectation value of $\mathcal{E}_{R_n(\delta)}(\cdot,cn)$.

\begin{prop}\label{expect}
Let $c>\frac{1}{2}$  and $\delta>0$ be fixed. Then the expectation value of $\mathcal{E}_{R_n(\delta)}(\cdot,cn)$ satisfies
\begin{align*}
\mathbb E\big(\mathcal{E}_{R_n(\delta)}(\cdot,cn)\big)=\frac{\delta^{1-2c}}{2c-1}
\end{align*}
for all $n\geq1$.
\end{prop}

\begin{Proof}
Using Siegel's formula \eqref{siegelform} and equations \eqref{unitvolume} and \eqref{rnd}, we find that
\begin{align*}
\mathbb E\big(\mathcal{E}_{R_n(\delta)}(\cdot,cn)\big)&=V_n^{-2c}\int_{\R^n}|\vecx|^{-2cn}I_{R_n(\delta)}(\vecx)\,d\vecx=\omega_nV_n^{-2c}\int_{R_n(\delta)}^{\infty}r^{n(1-2c)-1}\,dr\nonumber\\
&=\omega_nV_n^{-2c}\frac{R_n(\delta)^{n(1-2c)}}{n(2c-1)}=\frac{V_n^{1-2c}}{2c-1}\Big(\frac{\delta}{V_n}\Big)^{1-2c}=\frac{\delta^{1-2c}}{2c-1}
\end{align*}
for all $n\geq1$.
\end{Proof}

We continue by discussing the variance of $\mathcal{E}_{R_n(\delta)}(\cdot,cn)$.

\begin{prop}\label{var}
Let $c>\frac{1}{2}$  and $\delta>0$ be fixed. Then the variance of $\mathcal{E}_{R_n(\delta)}(\cdot,cn)$ satisfies
\begin{align*}
V\big(\mathcal E_{R_n(\delta)}(\cdot,cn)\big)\to\frac{2\delta^{1-4c}}{4c-1}
\end{align*}
as $n\to\infty$.
\end{prop}

\begin{Proof}
Following \cite{sst} we find that
\begin{align*}
&\mathbb E\big(\mathcal E_{R_n(\delta)}(\cdot,cn)^2\big)=\int_{X_n}\mathcal E_{R_n(\delta)}(L,cn)^2\,d\mu_n(L)\nonumber\\
&=V_n^{-4c}\int_{X_n}\sum_{\vecm_1\in L}|\vecm_1|^{-2cn}I_{R_n(\delta)}(\vecm_1)\sum_{\vecm_2\in L}|\vecm_2|^{-2cn}I_{R_n(\delta)}(\vecm_2)\,d\mu_n(L)\nonumber\\
&=V_n^{-4c}\underset{\substack{d_1\geq1\\d_2\geq1}}{\sum}\int_{X_n}\underset{\vecnu_1,\vecnu_2\in L}{{\sum}^*}|d_1\vecnu_1|^{-2cn}I_{R_n(\delta)}(d_1\vecnu_1)|d_2\vecnu_2|^{-2cn}I_{R_n(\delta)}(d_2\vecnu_2)\,d\mu_n(L),
\end{align*}
where $*$ denotes that the summation is over pairs of primitive vectors\footnote{A primitive lattice vector is a non-zero lattice vector which is not a positive integral multiple of another lattice vector.} in $L$. Applying a version of Rogers' formula adapted to the present situation (cf.\ \cite[Thm.\ 5]{rogers1}) we get
\begin{align*}
&\mathbb E\big(\mathcal{E}_{R_n(\delta)}(\cdot,cn)^2\big)\nonumber\\
&=V_n^{-4c}\underset{\substack{d_1\geq1\\d_2\geq1}}{\sum}\bigg(\frac{1}{\zeta(n)^2}\int_{\R^n}\int_{\R^n}|d_1\vecx_1|^{-2cn}I_{R_n(\delta)}(d_1\vecx_1)|d_2\vecx_2|^{-2cn}I_{R_n(\delta)}(d_2\vecx_2)\,d\vecx_1d\vecx_2\nonumber\\
&\hspace{60pt}+\frac{2}{\zeta(n)}\int_{\R^n}|d_1\vecx|^{-2cn}I_{R_n(\delta)}(d_1\vecx)|d_2\vecx|^{-2cn}I_{R_n(\delta)}(d_2\vecx)\,d\vecx\bigg)\nonumber\\
&=V_n^{-4c}\underset{\substack{d_1\geq1\\d_2\geq1}}{\sum}\frac{1}{\zeta(n)^2d_1^nd_2^n}\int_{\R^n}|\vecx_1|^{-2cn}I_{R_n(\delta)}(\vecx_1)\,d\vecx_1\int_{\R^n}|\vecx_2|^{-2cn}I_{R_n(\delta)}(\vecx_2)\,d\vecx_2\nonumber\\
&+\frac{2V_n^{-4c}}{\zeta(n)}\underset{\substack{d_1\geq1\\d_2\geq1}}{\sum}d_1^{-2cn}d_2^{-2cn}\int_{\R^n}|\vecx|^{-4cn}I_{R_n(\delta)}(d_1\vecx)I_{R_n(\delta)}(d_2\vecx)\,d\vecx\nonumber\\
&=\mathbb E\big(\mathcal E_{R_n(\delta)}(\cdot,cn)\big)^2+\frac{2V_n^{-4c}}{\zeta(n)}\underset{\substack{d_1\geq1\\d_2\geq1}}{\sum}d_1^{-2cn}d_2^{-2cn}\int_{\R^n}|\vecx|^{-4cn}I_{R_n(\delta)}(d_1\vecx)I_{R_n(\delta)}(d_2\vecx)\,d\vecx.
\end{align*}
Using equations \eqref{unitvolume} and \eqref{rnd} it follows that
\begin{align*}
V\big(\mathcal{E}_{R_n(\delta)}(\cdot,cn)\big)&=\frac{2V_n^{-4c}}{\zeta(n)}\underset{\substack{d_1\geq1\\d_2\geq1}}{\sum}d_1^{-2cn}d_2^{-2cn}\int_{\R^n}|\vecx|^{-4cn}I_{R_n(\delta)}(d_1\vecx)I_{R_n(\delta)}(d_2\vecx)\,d\vecx\nonumber\\
&=\frac{2V_n^{-4c}}{\zeta(n)}\underset{\substack{d_1\geq1\\d_2\geq1}}{\sum}d_1^{-2cn}d_2^{-2cn}\omega_n\int_{R_n(\delta)/\min(d_1,d_2)}^{\infty}r^{n(1-4c)-1}\,dr\nonumber\\
&=\frac{2V_n^{1-4c}}{\zeta(n)(4c-1)}\underset{\substack{d_1\geq1\\d_2\geq1}}{\sum}d_1^{-2cn}d_2^{-2cn}\Big(\frac{R_n(\delta)}{\min(d_1,d_2)}\Big)^{n(1-4c)}\nonumber\\
&=\frac{2\delta^{1-4c}}{\zeta(n)(4c-1)}\underset{\substack{d_1\geq1\\d_2\geq1}}{\sum}\min(d_1,d_2)^{n(2c-1)}\max(d_1,d_2)^{-2cn}.
\end{align*}
The sum above is rapidly decaying and the term corresponding to $d_1=d_2=1$ is exponentially larger than the remaining series. We conclude that
\begin{align*}
V\big(\mathcal E_{R_n(\delta)}(\cdot,cn)\big)\to\frac{2\delta^{1-4c}}{4c-1}
\end{align*}
as $n\to\infty$.
\end{Proof}

We now turn our attention to the higher order moments. Let $k\geq2$. To begin with we find that
\begin{align*}
&\mathbb E\big(\mathcal E_{R_n(\delta)}(\cdot,cn)^k\big)=\int_{X_n}\mathcal E_{R_n(\delta)}(L,cn)^k\,d\mu_n(L)\nonumber\\
&=V_n^{-2kc}\int_{X_n}\underset{\vecm_1,\ldots,\vecm_k\in L}{\sum}|\vecm_1|^{-2cn}I_{R_n(\delta)}(\vecm_1)\cdots|\vecm_k|^{-2cn}I_{R_n(\delta)}(\vecm_k)\,d\mu_n(L).
\end{align*}
Applying Rogers' formula \eqref{rogform} yields
\begin{align}\label{Ek}
&\mathbb E\big(\mathcal E_{R_n(\delta)}(\cdot,cn)^k\big)=\mathbb E\big(\mathcal E_{R_n(\delta)}(\cdot,cn)\big)^k+\underset{(\nu,\mu)}{\sum}\sum_{q=1}^{\infty}\sum_D\Big(\frac{e_1}{q}\cdots\frac{e_m}{q}\Big)^nI(D,n,c,\delta)
\end{align}
where
\begin{align*}
I(D,n,c,\delta)=V_n^{-2kc}\int_{\R^n}\cdots\int_{\R^n}\prod_{j=1}^k\bigg(\Big|\sum_{i=1}^m\frac{d_{ij}}{q}\vecx_i\Big|^{-2cn}I_{R_n(\delta)}\Big(\sum_{i=1}^m\frac{d_{ij}}{q}\vecx_i\Big)\bigg)\,d\vecx_1\ldots d\vecx_m.
\end{align*}
For some special $k$-admissible matrices $D$, the integral $I(D,n,c,\delta)$ is easy to determine.

\begin{prop}\label{propi}
Let $k\geq2$ and $1\leq m\leq k-1$ be given and let $(k_1,k_2,\ldots,k_m)$ be an ordered partition of $k$. Let $X_{k_1,k_2,\ldots,k_m}$ denote the set of $k$-admissible $m\times k$ matrices $D$, with elements $d_{ij}\in\{0,\pm1\}$, having exactly one non-zero entry in each column and $k_i$ non-zero entries in the $i$:th row. Let $c>\frac{1}{2}$  and $\delta>0$ be fixed. Then
\begin{align*}
 I(D,n,c,\delta)=\delta^{m-2kc}\prod_{i=1}^m\frac{1}{2k_ic-1}
\end{align*}
for all $D\in X_{k_1,k_2,\ldots,k_m}$ and all $n\geq1$.
\end{prop}

\begin{Proof}
For $D\in X_{k_1,k_2,\ldots,k_m}$ it is immediate to get
\begin{align*}
I(D,n,c,\delta)=V_n^{-2kc}\prod_{i=1}^m\int_{\R^n}|\vecx_i|^{-2k_icn}I_{R_n(\delta)}(\vecx_i)\,d\vecx_i.
\end{align*}
Using equations \eqref{unitvolume} and \eqref{rnd} it follows that
\begin{align*}
I(D,n,c,\delta)&=V_n^{-2kc}\prod_{i=1}^m\Big(\omega_n\int_{R_n(\delta)}^{\infty}r^{n-2k_icn-1}\,dr\Big)\\
&=V_n^{-2kc}\prod_{i=1}^m\Big(\omega_n\frac{R_n(\delta)^{n(1-2k_ic)}}{n(2k_ic-1)}\Big)\\
&=V_n^{m-2kc}R_n(\delta)^{n(m-2kc)}\prod_{i=1}^m\frac{1}{2k_ic-1}\\
&=\delta^{m-2kc}\prod_{i=1}^m\frac{1}{2k_ic-1}
\end{align*}
for all $n\geq1$.
\end{Proof}

\begin{prop}\label{propx}
Let $k\geq2$ and $1\leq m\leq k-1$ be given and let $(k_1,k_2,\ldots,k_m)$ be an ordered partition of $k$. Then
\begin{align*}
\#X_{k_1,k_2,\ldots,k_m}=2^{k-m}\prod_{i=1}^{m-1}\binom{k-\sum_{j=1}^{i-1}k_j-1}{k_i-1},
\end{align*}
where the empty product is interpreted as $1$ and the empty sum as $0$.
\end{prop}

\begin{Proof}
First we note that it is enough to count the number of relevant matrices having all non-zero entries equal to 1, since we get $\#X_{k_1,k_2,\ldots,k_m}$ as $2^{k-m}$ times this number.

The positions of the $k_1$ non-zero entries in the first row can be chosen in $\binom{k-1}{k_1-1}$ ways, since we must have $d_{1,1}=1$. Given such a configuration, the positions of the $k_2$ non-zero entries in the second row can be chosen in $\binom{k-k_1-1}{k_2-1}$ ways, since the first non-zero entry in the second row is determined by the configuration in the first row. Continuing the argument in the same way we get the desired result.
\end{Proof}

We now state our main result on the moments of $\mathcal{E}_{R_n(\delta)}(\cdot,cn)$.

\begin{thm}\label{higher}
Let $k\geq2$ and $1\leq m\leq k-1$. Let $c>\frac{1}{2}$  and $\delta>0$ be fixed. Then the $k$:th moment of $\mathcal E_{R_n(\delta)}(\cdot,cn)$ satisfies
\begin{align}\label{mainformula}
\mathbb E\big(\mathcal E_{R_n(\delta)}(\cdot,cn)^k\big)\to \sum_{m=1}^{k}\underset{k_1+\cdots+k_m=k}{\sum}2^{k-m}\prod_{i=1}^{m-1}\binom{k-\sum_{j=1}^{i-1}k_j-1}{k_i-1}\delta^{m-2kc}\prod_{i=1}^m\frac{1}{2k_ic-1}
\end{align}
as $n\to\infty$. In \eqref{mainformula} the inner sum is over all ordered partitions of $k$ in $m$ parts.
\end{thm}

We postpone the proof of Theorem \ref{higher} until Section \ref{heavy}.

\begin{remark}\label{higherremark}
The contribution to the right hand side of \eqref{mainformula} with $m=k$ and $k_1=k_2=\cdots=k_m=1$ comes from the term $\mathbb E\big(\mathcal E_{R_n(\delta)}(\cdot,cn)\big)^k$ in \eqref{Ek}. The remaining part of the right hand side of \eqref{mainformula} is the contribution from the matrices $D$ belonging to $X_{k_1,k_2,\ldots,k_m}$ for some ordered partition $(k_1,k_2,\ldots,k_m)$ of $k$.   
\end{remark}

\begin{remark}
By comparison with Proposition \ref{expect} we note that formula \eqref{mainformula} holds also with $k=1$. We further note the agreement of Theorem \ref{higher} with Proposition \ref{var} (using the appropriate combination of $k=1$ and $k=2$).
\end{remark}

\section{Proof of Theorem \ref{mainthm}}\label{mainsec}

Let $\frac{1}{2}<c_1<\cdots<c_m$ be fixed and set $\vecc=(c_1,\ldots,c_m)$. We want to prove that the distribution of the random vector 
\begin{align}\label{Evec}
\mathcal E_n(L,\vecc n):=\big(\mathcal E_n(L,c_1n),\ldots,\mathcal E_n(L,c_mn)\big) 
\end{align} 
converges to the distribution of 
\begin{align*}
T(\vecc)=\big(T(c_1),\ldots,T(c_m)\big)
\end{align*}
as $n\to\infty$. This will be done working with partial sums of $\mathcal E_n(L,c)$ and $T(c)$ respectively. For $k\geq1$ and $L\in X_n$ we let
\begin{align*}
\mathcal E_n(L,\vecc n,k):=\Big(2\sum_{j=1}^{k}\mathcal{V}_j(L)^{-2c_1},\ldots,2\sum_{j=1}^{k}\mathcal{V}_j(L)^{-2c_m}\Big) 
\end{align*}
and
\begin{align*}
T(\vecc,k):=\Big(2\sum_{j=1}^{k}T_j^{-2c_1},\ldots,2\sum_{j=1}^{k}T_j^{-2c_m}\Big).
\end{align*}

\begin{lem}\label{easylemma}
Let $k\geq1$ be fixed. Then $\mathcal E_n(\cdot,\vecc n,k)$ converges in distribution to $T(\vecc,k)$ as $n\to\infty$.
\end{lem}

\begin{Proof}
From Theorem \ref{thmjag} it follows that the random vector $\big(\mathcal{V}_1(\cdot),\ldots,\mathcal{V}_k(\cdot)\big)$ converges in distribution to the random vector $(T_1,\ldots,T_k)$ as $n\to\infty$. Since \mbox{$f:(\R^+)^k\to\R^m$} defined by  
\begin{align*}
f(x_1,\ldots,x_k)=\Big(2\sum_{j=1}^{k}x_j^{-2c_1},\ldots,2\sum_{j=1}^{k}x_j^{-2c_m}\Big) 
\end{align*}
is continuous the lemma follows from \cite[Thm.\ 2.7 or simply eq.\ (2.5)]{billconv}. 
\end{Proof}

We also note that by definition $\lim_{k\to\infty}\mathcal E_n(L,\vecc n,k)=\mathcal E_n(L,\vecc n)$ for each fixed \mbox{$L\in X_n$}. Since sure convergence implies convergence in distribution we get that $\mathcal E_n(\cdot,\vecc n,k)$ converges in distribution to $\mathcal E_n(\cdot,\vecc n)$ as $k\to\infty$ for each fixed $n$. Similarly we find that $T(\vecc,k)$ converges in distribution to $T(\vecc)$ as $k\to\infty$. In short it is this observation, but in the more precise form given in Proposition \ref{levy}, together with Lemma \ref{easylemma} that proves Theorem \ref{mainthm}.

We will during the remaining part of the proof of Theorem \ref{mainthm} change our point of view on convergence in distribution. We let $\mathcal P(\R^m)$ denote the set of Borel probability measures on $\R^m$. We recall that for $P,Q\in\mathcal P(\R^m)$ the Lévy-Prohorov distance $\pi(P,Q)$ between $P$ and $Q$ is defined as 
\begin{align}\label{LP}
 \pi(P,Q):=\inf\big\{\ve>0\mid P(A)\leq Q(A^{\ve})+\ve\,\text{ for all Borel sets $A\subseteq\R^m$}\big\},
\end{align}
where $A^{\ve}$ is the open $\ve$-neighbourhood of $A$ in $\R^m$ (cf.\ \cite{billconv}). Since $\R^m$ is separable, it is known that convergence in the metric $\pi$ is equivalent to weak convergence in $\mathcal P(\R^m)$. We prove the following:

\begin{prop}\label{levy}
Let $\vecc=(c_1,\ldots,c_m)$ be fixed as above. Let $\mu_{\mathcal E_n(\cdot,\vecc n,k)}$, $\mu_{\mathcal E_n(\cdot,\vecc n)}$, $\mu_{T(\vecc,k)}$ and $\mu_{T(\vecc)}$ be the distributions of the random vectors $\mathcal E_n(\cdot,\vecc n,k)$, $\mathcal E_n(\cdot,\vecc n)$, $T(\vecc,k)$ and $T(\vecc)$ respectively. Then for every $\ve>0$ there exists $K,N\in\Z^+$ such that 
\begin{align}\label{prop1}
\pi(\mu_{\mathcal E_n(\cdot,\vecc n,k)},\mu_{\mathcal E_n(\cdot,\vecc n)})\leq\ve 
\end{align}
and 
\begin{align}\label{prop2}
 \pi(\mu_{T(\vecc,k)},\mu_{T(\vecc)})\leq\ve
\end{align}
for all $k\geq K$ and all $n\geq N$.
\end{prop}

\begin{Proof}
Let $\ve>0$ be given. We note that for $Y>0$ we have 
\begin{align}\label{complex}
&\mathbb E\bigg(\Big|\Big(2\sum_{\mathcal{V}_j>Y}\mathcal{V}_j(\cdot)^{-2c_1},\ldots,2\sum_{\mathcal{V}_j>Y}\mathcal{V}_j(\cdot)^{-2c_m}\Big)\Big|^2\bigg)\nonumber\\
&=\int_{X_n}\Big(\mathcal E_{R_n(Y)}(L,c_1n)^2+\ldots+\mathcal E_{R_n(Y)}(L,c_mn)^2\Big)\,d\mu_n(L). 
\end{align}
Hence, by Theorem \ref{higher} with $k=2$, we find that
\begin{align}\label{secondmoment}
&\mathbb E\bigg(\Big|\Big(2\sum_{\mathcal{V}_j>Y}\mathcal{V}_j(\cdot)^{-2c_1},\ldots,2\sum_{\mathcal{V}_j>Y}\mathcal{V}_j(\cdot)^{-2c_m}\Big)\Big|^2\bigg)\to\sum_{i=1}^{m}\bigg(\Big(\frac{Y^{1-2c_i}}{2c_i-1}\Big)^2+2\frac{Y^{1-4c_i}}{4c_i-1}\bigg)\end{align}
as $n\to\infty$. We now fix $Y_0$ large enough to make the right hand side in \eqref{secondmoment}, with $Y=Y_0$, less than $\frac{1}{4}\ve^3$. Then it follows from \eqref{secondmoment} that there exists $N_{\ve}\in\Z^+$ such that 
\begin{align}\label{star}
&\mathbb E\bigg(\Big|\Big(2\sum_{\mathcal{V}_j>Y_0}\mathcal{V}_j(\cdot)^{-2c_1},\ldots,2\sum_{\mathcal{V}_j>Y_0}\mathcal{V}_j(\cdot)^{-2c_m}\Big)\Big|^2\bigg)<\frac{\ve^3}{2} 
\end{align}
for all $n\geq N_{\ve}$.

We next study $\text{Prob}_{\mu_n}\big\{\#\{\mathcal V_j(L)\leq Y_0\}>K\big\}$ with $K\in\Z^+$ and $n\geq N_{\ve}$. By \cite[Thm.\ 3]{rogers3} (or our Theorem \ref{thmjag}) and by possibly increasing $N_{\ve}$ we can choose $K_{\ve}\in\Z^+$, depending on $Y_0$ and $N_{\ve}$, such that \begin{align}\label{square} 
\text{Prob}_{\mu_n}\big\{\#\{\mathcal V_j(L)\leq Y_0\}>K_{\ve}\big\}<\frac{\ve}{2}.                                         \end{align}
For fixed $L\in X_n$ ($n\geq N_{\ve}$) either $\#\{\mathcal V_j(L)\leq Y_0\}>K_{\ve}$ 
or $\#\{\mathcal V_j(L)\leq Y_0\}\leq K_{\ve}$. We let 
\begin{align*}
\widetilde X_n:=\big\{L\in X_n\mid \#\{\mathcal V_j(L)\leq Y_0\}\leq K_{\ve}\big\}. 
\end{align*}
Using \eqref{star} we get, for $k\geq K_{\ve}$,
\begin{align}\label{starstar}
&\int_{\widetilde X_n}I\Big(\big|\mathcal E_n(\cdot,\vecc n)-\mathcal E_n(\cdot,\vecc n,k)\big|\geq\ve\Big)\,d\mu_n(L)\nonumber\\
&\leq\frac{1}{\ve^2}\int_{\widetilde X_n}\big|\mathcal E_n(\cdot,\vecc n)-\mathcal E_n(\cdot,\vecc n,k)\big|^2\,d\mu_n(L)\nonumber\\
&=\frac{1}{\ve^2}\int_{\widetilde X_n}\Big|\Big(2\sum_{j=k+1}^{\infty}\mathcal{V}_j(\cdot)^{-2c_1},\ldots,2\sum_{j=k+1}^{\infty}\mathcal{V}_j(\cdot)^{-2c_m}\Big)\Big|^2\,d\mu_n(L)\nonumber\\
&\leq\frac{1}{\ve^2}\mathbb E\bigg(\Big|\Big(2\sum_{\mathcal{V}_j>Y_0}\mathcal{V}_j(\cdot)^{-2c_1},\ldots,2\sum_{\mathcal{V}_j>Y_0}\mathcal{V}_j(\cdot)^{-2c_m}\Big)\Big|^2\bigg)<\frac{\ve}{2}.
\end{align}
For $n\geq N_{\ve}$ and $k\geq K_{\ve}$ it follows from \eqref{square} and \eqref{starstar} that for every Borel set $A\subseteq\R^m$ we have
\begin{align*}
\mu_{\mathcal E_n(\cdot,\vecc n,k)}(A)<\mu_{\mathcal E_n(\cdot,\vecc n)}(A^{\ve})+\frac{\ve}{2}+\frac{\ve}{2}=\mu_{\mathcal E_n(\cdot,\vecc n)}(A^{\ve})+\ve,
\end{align*}
which in turn implies that $\pi(\mu_{\mathcal E_n(\cdot,\vecc n,k)},\mu_{\mathcal E_n(\cdot,\vecc n)})\leq\ve$. Hence inequality \eqref{prop1} holds for all $k\geq K_{\ve}$ and $n\geq N_{\ve}$.

Finally we recall from the discussion following Lemma \ref{easylemma} that $T(\vecc,k)$ converges in distribution to $T(\vecc)$ as $k\to\infty$. Hence, by possibly increasing $K_{\ve}$, we find that also inequality  \eqref{prop2} holds. We conclude that the statement of the proposition holds with $K=K_{\ve}$ and $N=N_{\ve}$.
\end{Proof}

\begin{Proof1} 
Let $\ve>0$ be given and let $K_{\ve}$ and $N_{\ve}$ be as in the proof of Proposition \ref{levy}. It follows from Lemma \ref{easylemma} that there exists $N_0\in\Z^+$ such that $\pi(\mu_{\mathcal E_n(\cdot,\vecc n,K_{\ve})},\mu_{T(\vecc,K_{\ve})})<\ve$ for all $n\geq N_{0}$. This together with Proposition \ref{levy} yields
\begin{align*}
\pi(\mu_{\mathcal E_n(\cdot,\vecc n)},\mu_{T(\vecc)})&\leq\pi(\mu_{\mathcal E_n(\cdot,\vecc n)},\mu_{\mathcal E_n(\cdot,\vecc n,K_{\ve})})\\
&+
\pi(\mu_{\mathcal E_n(\cdot,\vecc n,K_{\ve})},\mu_{T(\vecc,K_{\ve})})+
\pi(\mu_{T(\vecc,K_{\ve})},\mu_{T(\vecc)})<3\ve 
\end{align*}
for all $n\geq\max(N_\ve,N_0)$. We conclude that $\mu_{\mathcal E_n(\cdot,\vecc n)}$ converges in the metric $\pi$ to $\mu_{T(\vecc)}$ as $n\to\infty$. From the discussion just above Proposition \ref{levy} we know that this is equivalent to that $\mathcal E_n(\cdot,\vecc n)$ converges in distribution to $T(\vecc)$ as $n\to\infty$.
\end{Proof1}

\begin{remark}
Theorem \ref{mainthm} holds also with $\vecc=(c_1,\ldots,c_m)\in\C^m$ satisfying $\Re c_i>\frac{1}{2}$, $1\leq i\leq m$. The proof is the same except that in the second line of \eqref{complex} we get a sum of the second \emph{absolute} moments of $\mathcal E_{R_n(Y)}(\cdot,c_in)$, $1\leq i\leq m$. The limits, as $n\to\infty$, of these absolute moments can be calculated in the same way as the limits in Theorem \ref{higher}; in particular the right hand side of \eqref{secondmoment} turns into 
\begin{align*}
\sum_{i=1}^{m}\bigg(\Big(\frac{Y^{1-2\Re c_i}}{|2c_i-1|}\Big)^2+2\frac{Y^{1-4\Re c_i}}{4\Re c_i-1}\bigg).
\end{align*}

In a related vein we mention that it is also possible to determine, for any given $c>\frac{1}{2}$, the limit as $n\to\infty$ of the probability of $E_n(L,s)$ having a (complex) zero $s$ with $\Re s>cn$. This result, which is joint work with A. Str\"ombergsson, will be presented elsewhere.
\end{remark}

\section{Proof of Theorem \ref{higher}}\label{heavy}

\subsection{Some first estimates}

In this section we follow Rogers (\cite[Sec.\ 9]{rogers2}). We begin by proving an upper bound for $I(D,n,c,\delta)$.

\begin{lem}\label{Iest}
Let $k\geq2$. Let $c>\frac{1}{2}$ and $\delta>0$ be fixed. Then, for any $k$-admissible $m\times k$ matrix $D$, we have
\begin{align*}
 I(D,n,c,\delta)\leq M(D)^{-n}\frac{\delta^{m-2kc}}{(2c-1)^m}
\end{align*}
where $M(D)$ is $q^{-m}$ times the largest value taken by any determinant of an $m\times m$-minor of $D$.
\end{lem}

\begin{Proof}
If $\lambda_1,\ldots,\lambda_m$ are the indices of any choice of $m$ linearly independent columns of $D$, then
\begin{align*}
\vecz_j=\sum_{i=1}^m\frac{d_{i\lambda_j}}{q}\vecx_i\hspace{15pt}(1\leq j\leq m)
\end{align*}
defines a linear change of variables with determinant of absolute value 
\begin{align*}
\mathcal D=\Big|\det\Big(\frac{d_{i\lambda_j}}{q}\Big)_{i,j=1}^m\Big|^n. 
\end{align*}
Using equations \eqref{unitvolume} and \eqref{rnd} and estimating the factors in $I(D,n,c,\delta)$ coming from columns with indices not in $\{\lambda_1,\ldots,\lambda_m\}$ with their respective supremum over $\R^n$, we get
\begin{align*}
I(D,n,c,\delta)&\leq  V_n^{-2kc}R_n(\delta)^{-2(k-m)cn}\int_{\R^n}\cdots\int_{\R^n}\prod_{j=1}^m\big(|\vecz_j|^{-2cn}I_{R_n(\delta)}(\vecz_j)\big)\,d\vecx_1\ldots d\vecx_m\\
&=\mathcal{D}^{-1}V_n^{-2kc}R_n(\delta)^{-2(k-m)cn}\bigg(\int_{\R^n}|\vecz|^{-2cn}I_{R_n(\delta)}(\vecz)\,d\vecz\bigg)^m\nonumber\\
&=\mathcal{D}^{-1}V_n^{-2kc}R_n(\delta)^{-2(k-m)cn}\bigg(\omega_n\int_{R_n(\delta)}^{\infty}r^{n-2cn-1}\,dr\bigg)^m\nonumber\\
&=\mathcal{D}^{-1}V_n^{m-2kc}\frac{R_n(\delta)^{mn-2kcn}}{(2c-1)^m}=\mathcal{D}^{-1}\frac{\delta^{m-2kc}}{(2c-1)^m}.
\end{align*}
Since this estimate holds for all changes of variables of the above type we arrive at the desired conclusion.
\end{Proof}

We now give a bound on the contribution from most of the terms in \eqref{Ek}.

\begin{prop}\label{intermed}
Let $c>\frac{1}{2}$ and $\delta>0$ be fixed. Let $k\geq2$ and assume that $n>\underset{1\leq m\leq k-1}{\max}\big(m(k-m)+1\big)$. 
Then
\begin{align*}
\underset{(\nu,\mu)}{\sum}\sum_{q=1}^{\infty}\underset{\substack{D}}{\sum}\Big(\frac{e_1}{q}\cdots\frac{e_m}{q}\Big)^nI(D,n,c,\delta)=\underset{(\nu,\mu)}{\sum}\underset{\substack{D\\q=1\\M(D)=1}}{\sum}I(D,n,c,\delta)+R(k),   
\end{align*}
where the remainder term satisfies 
\begin{align*}
0\leq R(k)\ll 2^{-n}.
\end{align*}
The implied constant depends on $c$, $\delta$ and $k$ but not on $n$. 
\end{prop}

\begin{Proof}
This is a straightforward adaptation of \cite[Sec.\ 9]{rogers2} to the present setting using inequality \eqref{element} and Lemma \ref{Iest}.
\end{Proof}

We note that every matrix $D$ with $q=1$ and $M(D)=1$ has all entries $d_{ij}\in\{0,\pm1\}$. In particular all the matrices in Proposition \ref{propi} are on this form. In Section \ref{finalest} we will discuss the contribution to \eqref{Ek} coming from matrices $D$ with $q=1$ and $M(D)=1$ and at least one column containing more than one non-zero entry. First we need to prove a series of integral estimates.

\subsection{Spherical symmetrization and integral estimates}\label{spherical}

In this section we let, for $c>\frac{1}{2}$ and $\delta>0$, 
\begin{align*}
f_{c,\delta}(\vecx)=|\vecx|^{-2cn}I_{R_n(\delta)}(\vecx).
\end{align*}
In order to use methods developed by Rogers (\cite{rogers3}, \cite{rogers4}) we first need to determine the function ${f_{c,\delta}}^*$, obtained from the function $f_{c,\delta}$ by spherical symmetrization. We let $\lambda$ denote the Lebesgue measure on $\R^n$. By the definition of ${f_{c,\delta}}^*$ (cf. \cite{rogers3}) we have
\begin{align*}
{f_{c,\delta}}^*(\vec0)=\underset{\vecx\in\R^n}{\sup}f_{c,\delta}(\vecx)=R_n(\delta)^{-2cn}
\end{align*}
and, for $\vecx\neq\vec0$,
\begin{align}\label{sym}
{f_{c,\delta}}^*(\vecx)&=\inf\big\{\rho>0\mid \lambda\big(\{\vecy : f_{c,\delta}(\vecy)>\rho\}\big)\leq\lambda\big(\{\vecy : |\vecy|\leq|\vecx|\}\big)\big\}\nonumber\\
&=\inf\Big\{\rho>0 \,\Big|\, \omega_n\int_{R_n(\delta)}^{\rho^{-\frac{1}{2cn}}}r^{n-1}\,dr\leq V_n|\vecx|^n\Big\}\nonumber\\
&=\inf\big\{\rho>0\mid \rho^{-\frac{1}{2c}}-R_n(\delta)^n\leq|\vecx|^n\big\}\nonumber\\
&=\inf\big\{\rho>0\mid \rho\geq(|\vecx|^n+R_n(\delta)^n)^{-2c}\big\}=\big(|\vecx|^n+R_n(\delta)^n\big)^{-2c}.
\end{align}
We note that the formula for ${f_{c,\delta}}^*$ in \eqref{sym} is valid also for $\vecx=\vec0$.

We next prove some technical integral estimates. First we estimate an integral involving ${f_{c,\delta}}^*$.

\begin{prop}\label{AS4}
Let $c>\frac{1}{2}$ and $\delta>0$ be fixed. Let $\ell\geq4$ and let $\ell_1,\ell_2,\ell_3,\ell_4$ be positive integers satisfying $\ell_1+\ell_2+\ell_3+\ell_4=\ell$. Then
\begin{align}\label{symint}
&V_n^{-2\ell c}\int_{\R^n}\int_{\R^n}{f_{c,\delta}}^*(\vecx_1)^{\ell_1}{f_{c,\delta}}^*(\vecx_2)^{\ell_2}{f_{c,\delta}}^*(\vecx_1+\vecx_2)^{\ell_3}{f_{c,\delta}}^*(\vecx_1-\vecx_2)^{\ell_4}\,d\vecx_1d\vecx_2\ll \mathcal C(n),
\end{align}
where $\mathcal{C}(n)$ decays exponentially with $n$. The implied constant depends on $\ell,c$ and $\delta$ but not on $n$.
\end{prop}

\begin{remark}
Our proof gives $\mathcal{C}(n)\ll\sqrt{n}\big(\frac{4}{5}\big)^{\frac{n}{2}}$, but we have not optimized to get the best possible constant. 
\end{remark}

\begin{Proof}
We call the left hand side in \eqref{symint} $I$. Changing to spherical coordinates and using the law of cosines we obtain
\begin{align*}
I&=\omega_n\omega_{n-1}V_n^{-2\ell c}\int_{0}^{\infty}\big(r_1^n+R_n(\delta)^n\big)^{-2\ell _1c}r_1^{n-1}\int_{0}^{\infty}\big(r_2^n+R_n(\delta)^n\big)^{-2\ell _2c}r_2^{n-1}\\
&\times\int_{0}^{\pi}\big((r_1^2+r_2^2+2r_1r_2\cos(\varphi))^{\frac{n}{2}}+R_n(\delta)^n\big)^{-2\ell _3c}\\
&\times\big((r_1^2+r_2^2-2r_1r_2\cos(\varphi))^{\frac{n}{2}}+R_n(\delta)^n\big)^{-2\ell _4c}\sin^{n-2}(\varphi)\,d\varphi dr_2dr_1\\
&=I(\ell _1,\ell _2,\ell _3,\ell _4)+I(\ell _1,\ell _2,\ell _4,\ell _3),
\end{align*}
where
\begin{align*}
&I(\ell _1,\ell _2,\ell _3,\ell _4)\\&=\omega_n\omega_{n-1}V_n^{-2\ell c}\int_{0}^{\infty}\big(r_1^n+R_n(\delta)^n\big)^{-2\ell _1c}r_1^{n-1}\int_{0}^{\infty}\big(r_2^n+R_n(\delta)^n\big)^{-2\ell _2c}r_2^{n-1}\\
&\times\int_{0}^{\frac{\pi}{2}}\big((r_1^2+r_2^2+2r_1r_2\cos(\varphi))^{\frac{n}{2}}+R_n(\delta)^n\big)^{-2\ell _3c}\\
&\times\big((r_1^2+r_2^2-2r_1r_2\cos(\varphi))^{\frac{n}{2}}+R_n(\delta)^n\big)^{-2\ell _4c}\sin^{n-2}(\varphi)\,d\varphi dr_2dr_1.
\end{align*}
By crude estimates of the last factors we get
\begin{align}\label{double}
&I(\ell _1,\ell _2,\ell _3,\ell _4)\nonumber\\
&\ll\omega_n\omega_{n-1}V_n^{-2\ell c}R_n(\delta)^{-2\ell _4cn}\int_{0}^{\infty}\big(r_1^n+R_n(\delta)^n\big)^{-2\ell _1c}r_1^{n-1}\int_{0}^{\infty}\big(r_2^n+R_n(\delta)^n\big)^{-2\ell _2c}r_2^{n-1}\nonumber\\
&\times\int_{0}^{\frac{\pi}{2}}\big((r_1^2+r_2^2+2r_1r_2\cos(\varphi))^{\frac{n}{2}}+R_n(\delta)^n\big)^{-2\ell _3c}\sin^{n-2}(\varphi)\,d\varphi dr_2dr_1\nonumber\\
&\ll\omega_n\omega_{n-1}V_n^{-2\ell c}R_n(\delta)^{-2\ell _4cn}\int_{0}^{\infty}\big(r_1^n+R_n(\delta)^n\big)^{-2\ell _1c}r_1^{n-1}\int_{r_1}^{\infty}\big(r_2^n+R_n(\delta)^n\big)^{-2\ell _2c}r_2^{n-1}\nonumber\\
&\times\big((r_1^2+r_2^2)^{\frac{n}{2}}+R_n(\delta)^n\big)^{-2\ell _3c}\,dr_2dr_1\nonumber\\
&+\omega_n\omega_{n-1}V_n^{-2\ell c}R_n(\delta)^{-2\ell _4cn}\int_{0}^{\infty}\big(r_1^n+R_n(\delta)^n\big)^{-2\ell _2c}r_1^{n-1}\int_{r_1}^{\infty}\big(r_2^n+R_n(\delta)^n\big)^{-2\ell _1c}r_2^{n-1}\nonumber\\
&\times\big((r_1^2+r_2^2)^{\frac{n}{2}}+R_n(\delta)^n\big)^{-2\ell _3c}\,dr_2dr_1.
\end{align}
We note that when $r_1\leq r_2\leq 2r_1$ we have
\begin{align}\label{Mineq}
r_1^2+r_2^2\geq\frac{5}{4}r_2^2.
\end{align}
Hence
\begin{align}\label{snart}
&\omega_n\omega_{n-1}V_n^{-2\ell c}R_n(\delta)^{-2\ell _4cn}\int_{0}^{\infty}\big(r_1^n+R_n(\delta)^n\big)^{-2\ell _1c}r_1^{n-1}\int_{r_1}^{\infty}\big(r_2^n+R_n(\delta)^n\big)^{-2\ell _2c}r_2^{n-1}\nonumber\\
&\times\big((r_1^2+r_2^2)^{\frac{n}{2}}+R_n(\delta)^n\big)^{-2\ell _3c}\,dr_2dr_1\nonumber\\
&\ll\omega_n\omega_{n-1}V_n^{-2\ell c}R_n(\delta)^{-2\ell _4cn}\int_{0}^{R_n(\delta)}\big(r_1^n+R_n(\delta)^n\big)^{-2\ell _1c}r_1^{n-1}\nonumber\\
&\times\int_{r_1}^{2r_1}\big(r_2^n+R_n(\delta)^n\big)^{-2\ell _2c}\big((\sfrac{5}{4})^{\frac{n}{2}}r_2^n+R_n(\delta)^n\big)^{-2\ell _3c}r_2^{n-1}\,dr_2dr_1\nonumber\\
&+\omega_n\omega_{n-1}V_n^{-2\ell c}R_n(\delta)^{-2\ell _4cn}\int_{0}^{R_n(\delta)}\big(r_1^n+R_n(\delta)^n\big)^{-2\ell _1c}r_1^{n-1}\nonumber\\
&\times\int_{2r_1}^{\infty}\big(r_2^n+R_n(\delta)^n\big)^{-2(\ell _2+\ell _3)c}r_2^{n-1}\,dr_2dr_1\nonumber\\
&+\omega_n\omega_{n-1}V_n^{-2\ell c}R_n(\delta)^{-2\ell _4cn}\int_{R_n(\delta)}^{\infty}\big(r_1^n+R_n(\delta)^n\big)^{-2\ell _1c}r_1^{n-1}\nonumber\\
&\times\int_{r_1}^{\infty}\big(r_2^n+R_n(\delta)^n\big)^{-2\ell _2c}r_2^{n-1}\big((r_1^2+r_2^2)^{\frac{n}{2}}+R_n(\delta)^n\big)^{-2\ell _3c}\,dr_2dr_1.
\end{align}
We call the three terms in the right hand side of \eqref{snart} $I_1$, $I_2$ and $I_3$ respectively.

We first bound $I_3$. Using \eqref{Mineq} once more we find that
\begin{align*}
I_3&\ll\omega_n\omega_{n-1}V_n^{-2\ell c}R_n(\delta)^{-2\ell _4cn}\int_{R_n(\delta)}^{\infty}r_1^{n-2\ell _1cn-1}\int_{r_1}^{\infty}r_2^{n-2\ell _2cn-1}\big(r_1^2+r_2^2\big)^{-\ell _3cn}\,dr_2dr_1\nonumber\\
&\ll\omega_n\omega_{n-1}V_n^{-2\ell c}R_n(\delta)^{-2\ell _4cn}\int_{R_n(\delta)}^{\infty}r_1^{n-2\ell _1cn-1}\int_{r_1}^{2r_1}r_2^{n-2\ell _2cn-1}\big(\sfrac{5}{4}r_2^2\big)^{-\ell _3cn}\,dr_2dr_1\nonumber\\
&+\omega_n\omega_{n-1}V_n^{-2\ell c}R_n(\delta)^{-2\ell _4cn}\int_{R_n(\delta)}^{\infty}r_1^{n-2\ell _1cn-1}\int_{2r_1}^{\infty}r_2^{n-2(\ell _2+\ell_3)cn-1}\,dr_2dr_1\nonumber\\
&\ll\Big(\big(\sfrac{4}{5}\big)^{\ell _3cn}+2^{(1-2(\ell _2+\ell _3)c)n}\Big)\omega_{n-1}V_n^{1-2\ell c}R_n(\delta)^{-2\ell _4cn}\int_{R_n(\delta)}^{\infty}r_1^{2n-2(\ell _1+\ell _2+\ell _3)cn-1}\,dr_1\nonumber\\
&\ll\Big(\big(\sfrac{4}{5}\big)^{\ell _3cn}+2^{(1-2(\ell _2+\ell _3)c)n}\Big)\frac{\omega_{n-1}V_n^{1-2\ell c}R_n(\delta)^{2n-2\ell cn}}{n}.
\end{align*}
It follows from Stirling's formula that 
\begin{align*}
\omega_n\sim\Big(\frac{2\pi e}{n}\Big)^{\frac{n}{2}}\sqrt{\frac{n}{\pi}} \hspace{15pt}\text{as $n \to\infty$.}
\end{align*}
Hence 
\begin{align}\label{sterling}
\omega_{n-1}\ll n^{\frac{3}{2}}V_n, 
\end{align}
where the implied constant does not depend on $n$. Using \eqref{rnd} and \eqref{sterling} we conclude that
\begin{align}\label{I3}
I_3&\ll\sqrt{n}\Big(\big(\sfrac{4}{5}\big)^{\ell _3cn}+2^{(1-2(\ell _2+\ell _3)c)n}\Big)V_n^{2-2\ell c}R_n(\delta)^{2n-2\ell cn}\nonumber\\
&\ll\sqrt{n}\Big(\big(\sfrac{4}{5}\big)^{\ell _3cn}+2^{(1-2(\ell _2+\ell _3)c)n}\Big).
\end{align}

We next bound $I_2$. Changing variables we find that
\begin{align*}
I_2&=\frac{\omega_n\omega_{n-1}}{n}V_n^{-2\ell c}R_n(\delta)^{-2\ell _4cn}\int_{0}^{R_n(\delta)}\big(r_1^n+R_n(\delta)^n\big)^{-2\ell _1c}r_1^{n-1}\\
&\times\int_{(2r_1)^n}^{\infty}\big(v+R_n(\delta)^n\big)^{-2(\ell _2+\ell _3)c}\,dvdr_1\\
&\ll\omega_{n-1}V_n^{1-2\ell c}R_n(\delta)^{-2\ell _4cn}\int_{0}^{R_n(\delta)}\big(r_1^n+R_n(\delta)^n\big)^{-2\ell _1c}\\
&\times\big((2r_1)^n+R_n(\delta)^n\big)^{1-2(\ell _2+\ell _3)c}r_1^{n-1}\,dr_1\\
&\ll\frac{\omega_{n-1}V_n^{1-2\ell c}R_n(\delta)^{-2(\ell _1+\ell _4)cn}}{n}\int_{0}^{R_n(\delta)^n}\big(2^nu+R_n(\delta)^n\big)^{1-2(\ell _2+\ell _3)c}\,du\\
&\ll2^{-n}\frac{\omega_{n-1}V_n^{1-2\ell c}R_n(\delta)^{2n-2\ell cn}}{n}.
\end{align*}
Using \eqref{rnd} and \eqref{sterling} we obtain
\begin{align}\label{N12}
I_2\ll\sqrt{n}2^{-n}.
\end{align}
In a similar way we estimate $I_1$:
\begin{align*}
I_1&\ll\omega_n\omega_{n-1}V_n^{-2\ell c}R_n(\delta)^{-2\ell _4cn}\int_{0}^{R_n(\delta)}\big(r_1^n+R_n(\delta)^n\big)^{-2(\ell _1+\ell _2)c}r_1^{n-1}\\
&\times\int_{r_1}^{2r_1}\big((\sfrac{5}{4})^{\frac{n}{2}}r_2^n+R_n(\delta)^n\big)^{-2\ell _3c}r_2^{n-1}\,dr_2dr_1\\
&=\frac{\omega_n\omega_{n-1}}{n}V_n^{-2\ell c}R_n(\delta)^{-2\ell _4cn}\int_{0}^{R_n(\delta)}\big(r_1^n+R_n(\delta)^n\big)^{-2(\ell _1+\ell _2)c}r_1^{n-1}\\
&\times\int_{r_1^n}^{(2r_1)^n}\big((\sfrac{5}{4})^{\frac{n}{2}}v+R_n(\delta)^n\big)^{-2\ell _3c}\,dvdr_1\\
&\ll\big(\sfrac{4}{5}\big)^{\frac{n}{2}}\omega_{n-1}V_n^{1-2\ell c}R_n(\delta)^{-2\ell _4cn}\\
&\times\int_{0}^{R_n(\delta)}\big(r_1^n+R_n(\delta)^n\big)^{-2(\ell _1+\ell _2)c}\big((\sfrac{5}{4})^{\frac{n}{2}}r_1^n+R_n(\delta)^n\big)^{1-2\ell _3c}r_1^{n-1}\,dr_1\\
&\ll\big(\sfrac{4}{5}\big)^{\frac{n}{2}}\frac{\omega_{n-1}V_n^{1-2\ell c}R_n(\delta)^{-2\ell _4cn}}{n}\int_{0}^{R_n(\delta)^n}\big(u+R_n(\delta)^n\big)^{1-2(\ell _1+\ell _2+\ell _3)c}\,du\\
&\ll\big(\sfrac{4}{5}\big)^{\frac{n}{2}}\frac{\omega_{n-1}V_n^{1-2\ell c}R_n(\delta)^{2n-2\ell cn}}{n}.
\end{align*}
Now \eqref{rnd} and \eqref{sterling} yield
\begin{align}\label{N11}
I_1\ll\sqrt{n}\big(\sfrac{4}{5}\big)^{\frac{n}{2}}.
\end{align}

Combining \eqref{I3}, \eqref{N12} and \eqref{N11} with the corresponding estimates for the last integral in \eqref{double} we get
\begin{align*}
 I(\ell _1,\ell _2,\ell _3,\ell _4)&\ll\sqrt{n}\Big(\big(\sfrac{4}{5}\big)^{\frac{n}{2}}+2^{-n}+\big(\sfrac{4}{5}\big)^{cn}+2^{(1-4c)n}\Big).
\end{align*}
Finally we obtain
\begin{align*}
I\ll\sqrt{n}\big(\sfrac{4}{5}\big)^{\frac{n}{2}},
\end{align*}
which proves the proposition.
\end{Proof}

In the spirit of Rogers (\cite{rogers3}) we can now prove the following theorem.

\begin{thm}\label{bra2}
Let $c>\frac{1}{2}$ and $\delta>0$ be fixed. Let $\ell\geq4$ and let $\ell_1,\ell_2,\ell_3,\ell_4$ be positive integers satisfying $\ell_1+\ell_2+\ell_3+\ell_4=\ell$. Let further $\vecy_i\in\R^n$ and $\ve_i\in\{1,-1\}$, $1\leq i\leq \ell-2$, be fixed. Then 
\begin{align*}
I_{\vecy}&:=V_n^{-2\ell c}\int_{\R^n}\int_{\R^n}f_{c,\delta}(\vecx_1)\Big(\prod_{i=1}^{\ell _1-1}f_{c,\delta}(\ve_i\vecx_1+\vecy_i)\Big)f_{c,\delta}(\vecx_2)\Big(\prod_{i=\ell _1}^{\ell _1+\ell _2-2}f_{c,\delta}(\ve_i\vecx_2+\vecy_i)\Big)\\
&\times\Big(\prod_{i=\ell _1+\ell _2-1}^{\ell _1+\ell _2+\ell _3-2}f_{c,\delta}\big(\ve_i(\vecx_1+\vecx_2)+\vecy_i\big)\Big)\Big(\prod_{i=\ell _1+\ell _2+\ell _3-1}^{\ell-2}f_{c,\delta}\big(\ve_i(\vecx_1-\vecx_2)+\vecy_i\big)\Big)\,d\vecx_1d\vecx_2\\
&\ll\mathfrak{C}(n),
\end{align*}
where $\mathfrak{C}(n)$ decays exponentially with $n$. The implied constant depends on $\ell,c$ and $\delta$ but not on $n,\ve_i$ and $\vecy_i$, $1\leq i\leq \ell-2$.
\end{thm}

\begin{Proof}
 Since ${f_{c,\delta}}^*(\vecx)$ is the spherical symmetrization of $f_{c,\delta}(\ve\vecx+\vecy)$ for any (fixed) $\vecy\in\R^n$ and $\ve\in\{1,-1\}$ it follows from \cite[Thm.\ 1]{rogers4} that
\begin{align*}
I_{\vecy}\leq V_n^{-2\ell c}\int_{\R^n}\int_{\R^n}{f_{c,\delta}}^*(\vecx_1)^{\ell _1}{f_{c,\delta}}^*(\vecx_2)^{\ell _2}{f_{c,\delta}}^*(\vecx_1+\vecx_2)^{\ell _3}{f_{c,\delta}}^*(\vecx_1-\vecx_2)^{\ell _4}\,d\vecx_1d\vecx_2.
\end{align*}
Hence the theorem follows from Proposition \ref{AS4}.
\end{Proof}

We continue with estimates similar to those in Proposition \ref{AS4} and Theorem \ref{bra2} for the case when $\ell_4=0$ (or $\ell_3=0$). The proofs are similar to the ones above and for this reason some parts will only be sketched.

\begin{prop}\label{AS2}
Let $c>\frac{1}{2}$ and $\delta>0$ be fixed. Let $\ell\geq3$ and let $\ell _1,\ell _2,\ell _3$ be positive integers satisfying $\ell _1+\ell _2+\ell _3=\ell$. Then
\begin{align*}
J:=V_n^{-2\ell c}\int_{\R^n}\int_{\R^n}{f_{c,\delta}}^*(\vecx_1)^{\ell_1}{f_{c,\delta}}^*(\vecx_2)^{\ell_2}{f_{c,\delta}}^*(\vecx_1+\vecx_2)^{\ell_3}\,d\vecx_1d\vecx_2
\ll \mathcal D(n),
\end{align*}
where $\mathcal{D}(n)$ decays exponentially with $n$. The implied constant depends on $\ell,c$ and $\delta$ but not on $n$.
 \end{prop}

\begin{Proof}
By changing variables we note that in order to estimate the size of $J$ it is enough to study the corresponding integral, $\widetilde J$, over the region where $|\vecx_1|,|\vecx_2|\leq|\vecx_1+\vecx_2|$. More precisely we have
\begin{align*}
 J\ll\widetilde J(\ell_1,\ell_2,\ell_3)+\widetilde J(\ell_2,\ell_3,\ell_1)+\widetilde J(\ell_3,\ell_1,\ell_2),
\end{align*}
where
\begin{align*}
\widetilde J(\ell_1,\ell_2,\ell_3)=V_n^{-2\ell c}\underset{|\vecx_1|,|\vecx_2|\leq|\vecx_1+\vecx_2|}{\int\int}{f_{c,\delta}}^*(\vecx_1)^{\ell_1}{f_{c,\delta}}^*(\vecx_2)^{\ell_2}{f_{c,\delta}}^*(\vecx_1+\vecx_2)^{\ell_3}\,d\vecx_1d\vecx_2.
\end{align*}

Passing to spherical coordinates we get
\begin{align*}
&\widetilde{ J}(\ell_1,\ell_2,\ell_3)\ll\omega_n\omega_{n-1}V_n^{-2\ell c}\int_{0}^{\infty}\big(r_1^n+R_n(\delta)^n\big)^{-2\ell _1c}r_1^{n-1}\int_{0}^{\infty}\big(r_2^n+R_n(\delta)^n\big)^{-2\ell _2c}r_2^{n-1}\\
&\times\int_{0}^{\theta}\Big(\big(\max\big(r_1,r_2,(r_1^2+r_2^2+2r_1r_2\cos(\varphi))^{\frac{1}{2}}\big)\big)^n+R_n(\delta)^n\Big)^{-2\ell _3c}\sin^{n-2}(\varphi)\,d\varphi dr_2dr_1\nonumber\\
&+\omega_n\omega_{n-1}V_n^{-2\ell c}\int_{0}^{\infty}\big(r_1^n+R_n(\delta)^n\big)^{-2\ell _1c}r_1^{n-1}\int_{0}^{\infty}\big(r_2^n+R_n(\delta)^n\big)^{-2\ell _2c}r_2^{n-1}\\
&\times\int_{\theta}^{\frac{2\pi}{3}}\big(\max(r_1,r_2)^n+R_n(\delta)^n\big)^{-2\ell _3c}\sin^{n-2}(\varphi)\,d\varphi dr_2dr_1,
\end{align*}
where $\theta=\arccos(-\sfrac{3}{16})$. We call the integrals above $\widetilde{J}_1$ and $\widetilde{J}_2$ respectively.

We split $\widetilde{J}_1$ as
\begin{align}\label{J1}
&\widetilde{J}_1=\omega_n\omega_{n-1}V_n^{-2\ell c}\int_{0}^{\infty}\big(r_1^n+R_n(\delta)^n\big)^{-2\ell _1c}r_1^{n-1}\int_{r_1}^{\infty}\big(r_2^n+R_n(\delta)^n\big)^{-2\ell _2c}r_2^{n-1}\nonumber\\
&\times\int_{0}^{\theta}\Big(\big(\max\big(r_1,r_2,(r_1^2+r_2^2+2r_1r_2\cos(\varphi))^{\frac{1}{2}}\big)\big)^n+R_n(\delta)^n\Big)^{-2\ell _3c}\sin^{n-2}(\varphi)\,d\varphi dr_2dr_1\nonumber\\
&+\omega_n\omega_{n-1}V_n^{-2\ell c}\int_{0}^{\infty}\big(r_1^n+R_n(\delta)^n\big)^{-2\ell _2c}r_1^{n-1}\int_{r_1}^{\infty}\big(r_2^n+R_n(\delta)^n\big)^{-2\ell _1c}r_2^{n-1}\nonumber\\
&\times\int_{0}^{\theta}\Big(\big(\max\big(r_1,r_2,(r_1^2+r_2^2+2r_1r_2\cos(\varphi))^{\frac{1}{2}}\big)\big)^n+R_n(\delta)^n\Big)^{-2\ell _3c}\sin^{n-2}(\varphi)\,d\varphi dr_2dr_1.
\end{align}
When $0\leq\varphi\leq\theta$ we have
\begin{align*}
 r_1^2+r_2^2+2r_1r_2\cos(\varphi)\geq r_1^2+r_2^2-\frac{3r_1r_2}{8}.
\end{align*}
In particular, if also $r_1\leq r_2\leq 2r_1$, we get
\begin{align}\label{ineq2}
 r_1^2+r_2^2+2r_1r_2\cos(\varphi)\geq \frac{r_1^2}{4}+r_2^2\geq\frac{17}{16}r_2^2.
\end{align}
Hence
\begin{align*}
&\omega_n\omega_{n-1}V_n^{-2\ell c}\int_{0}^{\infty}\big(r_1^n+R_n(\delta)^n\big)^{-2\ell _1c}r_1^{n-1}\int_{r_1}^{\infty}\big(r_2^n+R_n(\delta)^n\big)^{-2\ell _2c}r_2^{n-1}\\
&\times\int_{0}^{\theta}\Big(\big(\max\big(r_1,r_2,(r_1^2+r_2^2+2r_1r_2\cos(\varphi))^{\frac{1}{2}}\big)\big)^n+R_n(\delta)^n\Big)^{-2\ell _3c}\sin^{n-2}(\varphi)\,d\varphi dr_2dr_1\\
&\ll\omega_n\omega_{n-1}V_n^{-2\ell c}\int_{0}^{R_n(\delta)}\big(r_1^n+R_n(\delta)^n\big)^{-2\ell_1c}r_1^{n-1}\\
&\times\int_{r_1}^{2r_1}\big(r_2^n+R_n(\delta)^n\big)^{-2\ell _2c}\big((\sfrac{17}{16})^{\frac{n}{2}}r_2^n+R_n(\delta)^n\big)^{-2\ell _3c}r_2^{n-1}\,dr_2dr_1\\
&+\omega_n\omega_{n-1}V_n^{-2\ell c}\int_{0}^{R_n(\delta)}\big(r_1^n+R_n(\delta)^n\big)^{-2\ell _1c}r_1^{n-1}\int_{2r_1}^{\infty}\big(r_2^n+R_n(\delta)^n\big)^{-2(\ell _2+\ell _3)c}r_2^{n-1}\,dr_2dr_1\\
&+\omega_n\omega_{n-1}V_n^{-2\ell c}\int_{R_n(\delta)}^{\infty}\big(r_1^n+R_n(\delta)^n\big)^{-2\ell _1c}r_1^{n-1}\int_{r_1}^{\infty}\big(r_2^n+R_n(\delta)^n\big)^{-2\ell _2c}r_2^{n-1}\\
&\times\int_{0}^{\theta}\Big(\big(\max\big(r_1,r_2,(r_1^2+r_2^2+2r_1r_2\cos(\varphi))^{\frac{1}{2}}\big)\big)^n+R_n(\delta)^n\Big)^{-2\ell _3c}\sin^{n-2}(\varphi)\,d\varphi dr_2dr_1.
\end{align*}
We call the three integrals on the right hand side $\widetilde{J}_{1,1}$, $\widetilde{J}_{1,2}$ and $\widetilde{J}_{1,3}$ respectively. 

Using \eqref{rnd}, \eqref{sterling} and \eqref{ineq2} we find that 
\begin{align}\label{AJ13}
\widetilde{J}_{1,3}&\ll\omega_n\omega_{n-1}V_n^{-2\ell c}\int_{R_n(\delta)}^{\infty}r_1^{n-2\ell _1cn-1}\int_{r_1}^{\infty}r_2^{n-2\ell _2cn-1}\nonumber\\
&\times\int_{0}^{\theta}\max\big(r_1,r_2,(r_1^2+r_2^2+2r_1r_2\cos(\varphi))^{\frac{1}{2}}\big)^{-2\ell _3cn}\sin^{n-2}(\varphi)\,d\varphi dr_2dr_1\nonumber\\
&\ll\omega_n\omega_{n-1}V_n^{-2\ell c}\int_{R_n(\delta)}^{\infty}r_1^{n-2\ell_1cn-1}\int_{r_1}^{2r_1}r_2^{n-2\ell_2cn-1}\big(\sfrac{17}{16}r_2^2\big)^{-\ell _3cn}\,dr_2dr_1\nonumber\\
&+\omega_n\omega_{n-1}V_n^{-2\ell c}\int_{R_n(\delta)}^{\infty}r_1^{n-2\ell _1cn-1}\int_{2r_1}^{\infty}r_2^{n-2(\ell _2+\ell_3)cn-1}\,dr_2dr_1\nonumber\\
&\ll\sqrt{n}\Big(\big(\sfrac{16}{17}\big)^{\ell _3cn}+2^{(1-2(\ell_2+\ell_3)c)n}\Big).
\end{align}
In the same way as we estimated $I_2$ in the proof of Proposition \ref{AS4} we find that
\begin{align}\label{AJ12}
\widetilde{J}_{1,2}&\ll\sqrt{n}2^{-n},
\end{align}
and similarly we find that
\begin{align}\label{AJ11}
 \widetilde{J}_{1,1}&\ll\omega_n\omega_{n-1}V_n^{-2\ell c}\int_{0}^{R_n(\delta)}\big(r_1^n+R_n(\delta)^n\big)^{-2(\ell _1+\ell _2)c}r_1^{n-1}\nonumber\\
&\times\int_{r_1}^{2r_1}\big((\sfrac{17}{16})^{\frac{n}{2}}r_2^n+R_n(\delta)^n\big)^{-2\ell _3c}r_2^{n-1}\,dr_2dr_1
\ll\sqrt{n}\big(\sfrac{16}{17}\big)^{\frac{n}{2}}.
\end{align}
Combining \eqref{AJ13}, \eqref{AJ12} and \eqref{AJ11} with the corresponding estimates for the last integral in \eqref{J1} we get
\begin{align*}
 \widetilde{J}_1\ll\sqrt{n}\Big(\big(\sfrac{16}{17}\big)^{\frac{n}{2}}+2^{-n}+\big(\sfrac{16}{17}\big)^{cn}+2^{(1-4c)n}\Big).
\end{align*}

To finish the proof we estimate $\widetilde{J}_2$. Using \eqref{rnd} and \eqref{sterling} we find that
\begin{align*}
 \widetilde{J}_2&
\ll\sin^{n-2}(\theta)\omega_n\omega_{n-1}V_n^{-2\ell c}\int_{0}^{\infty}\big(r_1^n+R_n(\delta)^n\big)^{-2\ell _1c}r_1^{n-1}\\
&\times\int_{r_1}^{\infty}\big(r_2^n+R_n(\delta)^n\big)^{-2(\ell _2+\ell _3)c}r_2^{n-1}\,dr_2dr_1\\
&+\sin^{n-2}(\theta)\omega_n\omega_{n-1}V_n^{-2\ell c}\int_{0}^{\infty}\big(r_1^n+R_n(\delta)^n\big)^{-2\ell _2c}r_1^{n-1}\\
&\times\int_{r_1}^{\infty}\big(r_2^n+R_n(\delta)^n\big)^{-2(\ell _1+\ell _3)c}r_2^{n-1}\,dr_2dr_1\\
&\ll\sqrt{n}\sin^{n-2}(\theta).
\end{align*}
Hence
\begin{align*}
\widetilde{J}(\ell_1,\ell_2,\ell_3)&\ll\sqrt{n}\Big(\big(\sfrac{16}{17}\big)^{\frac{n}{2}}+2^{-n}+\big(\sfrac{16}{17}\big)^{cn}+2^{(1-4c)n}+\sin^{n-2}(\theta)\Big)
\end{align*}
and we conclude that
\begin{align*}
J&\ll\sqrt{n}\sin^{n-2}(\theta).
\end{align*}
This finishes the proof.
\end{Proof}

Proposition \ref{AS2} implies the following theorem.

\begin{thm}\label{bra4}
Let $c>\frac{1}{2}$ and $\delta>0$ be fixed. Let $\ell\geq3$ and let $\ell _1,\ell _2,\ell _3$ be positive integers satisfying $\ell _1+\ell _2+\ell _3=\ell$. Let further $\vecy_i\in\R^n$ and $\ve_i\in\{1,-1\}$, $1\leq i\leq \ell-2$, be fixed. Then
\begin{align*}
&V_n^{-2\ell c}\int_{\R^n}\int_{\R^n}f_{c,\delta}(\vecx_1)\Big(\prod_{i=1}^{\ell _1-1}f_{c,\delta}(\ve_i\vecx_1+\vecy_i)\Big)f_{c,\delta}(\vecx_2)\Big(\prod_{i=\ell _1}^{\ell _1+\ell _2-2}f_{c,\delta}(\ve_i\vecx_2+\vecy_i)\Big)\\
&\times\Big(\prod_{i=\ell _1+\ell _2-1}^{\ell-2}f_{c,\delta}\big(\ve_i(\vecx_1+\vecx_2)+\vecy_i\big)\Big)\,d\vecx_1d\vecx_2\ll \mathfrak{D}(n),
\end{align*}
where $\mathfrak{D}(n)$ decays exponentially with $n$. The implied constant depends on $\ell,c$ and $\delta$ but not on $n,\ve_i$ and $\vecy_i$, $1\leq i\leq \ell-2$.
\end{thm}

\subsection{The final estimate}\label{finalest}

In this section we will use the estimates from the previous section to obtain bounds on the contribution to \eqref{Ek} from the main terms in Proposition \ref{intermed} which are not among the terms treated in Proposition \ref{propi}.

\begin{prop}\label{last}
Let $c>\frac{1}{2}$ and $\delta>0$ be fixed. Let $k\geq3$ and fix $m$ satisfying $2\leq m\leq k-1$. Let $(\nu,\mu)$ be a division of the numbers $1,\ldots,k$ satisfying \eqref{division} with our $m$. Let $D$ be a $(\nu,\mu)$-admissible matrix with $q=1$, $M(D)=1$ and at least one column containing more than one non-zero entry. Then
\begin{align*}
I(D,n,c,\delta)\ll\mathfrak E(n),
\end{align*}
where $\mathfrak E(n)$ decays exponentially with $n$. The implied constant depends on $k,c$ and $\delta$ but not on $n$.
\end{prop}

\begin{Proof}
We let $1\leq\ell\leq k-m$ be such that the leftmost column with more than one non-zero entry is $\mu_{\ell}$ and we choose $\lambda_1$ and $\lambda_2$ to be minimal with the property that $\lambda_1<\lambda_2$ and $d_{\lambda_1\mu_{\ell}}\neq0\neq d_{\lambda_2\mu_{\ell}}$. We furthermore define the sets
\begin{align*}
&M=\{1,2,\ldots,m\}\setminus\{\lambda_1,\lambda_2\}\\
&N=\{j\mid d_{\lambda_1\mu_j}\neq0\,\text{ or }\, d_{\lambda_2\mu_j}\neq0\}\\
&P=\{1,\ldots,k-m\}\setminus N. 
\end{align*}
We recall that the assumptions on $D$ implies that $d_{ij}\in\{0,\pm1\}$ and write $I(D,n,c,\delta)$ as an iterated integral:
\begin{align}\label{inner}
&I(D,n,c,\delta)=V_n^{-2(k-\#N-2)c}\int_{\R^n}\cdots\int_{\R^n}\prod_{j\in M}f_{c,\delta}(\vecx_j)\prod_{j\in P}f_{c,\delta}\Big(\sum_{i=1}^md_{i\mu_j}\vecx_i\Big)\nonumber\\
&\times\bigg(V_n^{-2(\#N+2)c}\int_{\R^n}\int_{\R^n}f_{c,\delta}(\vecx_{\lambda_1})f_{c,\delta}(\vecx_{\lambda_2})\prod_{j\in N}f_{c,\delta}\Big(\sum_{i=1}^md_{i\mu_j}\vecx_i\Big)\,d\vecx_{\lambda_1}d\vecx_{\lambda_2}\bigg)\,\prod_{j\in M}d\vecx_j.
\end{align}
The inner integral in \eqref{inner} is either of the form in Theorem \ref{bra2} or of the form in Theorem \ref{bra4}. Hence
\begin{align}\label{outer}
&I(D,n,c,\delta)\nonumber\\
&\ll\mathfrak F(n)V_n^{-2(k-\#N-2)c}\int_{\R^n}\cdots\int_{\R^n}\prod_{j\in M}f_{c,\delta}(\vecx_j)\prod_{j\in P}f_{c,\delta}\Big(\sum_{i=1}^md_{i\mu_j}\vecx_i\Big)\,\prod_{j\in M}d\vecx_j,
\end{align}
where $\mathfrak F(n)$ decays exponentially with $n$ and the implied constant depends on $k,c$ and $\delta$ but not on $n$.

If the integral in \eqref{outer} is on the form in Proposition \ref{propi}, the bound \eqref{outer} implies the desired estimate. If, on the other hand, this is not the case we reach the same conclusion by induction.
\end{Proof}

Since there are only finitely many matrices of the form in Proposition \ref{last} (the number depends on $k$ and $m$ but not on $n$) this finishes the proof of Theorem \ref{higher}.

\section{The joint moments of $\mathcal E_{R_n(\delta)}(\cdot,\vecc n)$}\label{joint}

Let $m\geq1$ and fix $\vecc=(c_1,\ldots,c_m)$ satisfying $\frac{1}{2}<c_1< c_2<\cdots<c_m$. In this section we will outline what can be said about the joint moments of 
\begin{align}\label{Edvec}
\mathcal E_{R_n(\delta)}(L,\vecc n):=\big(\mathcal E_{R_n(\delta)}(L,c_1n),\ldots,\mathcal E_{R_n(\delta)}(L,c_mn)\big).  
\end{align}
The presentation here is parallel to the presentation in Section \ref{moments}.

Let $\kappa\geq2$  and fix $\vecgam=(\gamma_1,\ldots,\gamma_{\kappa})$ with $\gamma_j\in\{c_1,\ldots,c_m\}$, $1\leq j\leq\kappa$, satisfying $\gamma_1\leq \gamma_2\leq\cdots\leq \gamma_{\kappa}$. Applying Rogers' formula yields
\begin{align}\label{mix1}
&\mathbb E\Big(\prod_{j=1}^{\kappa}\mathcal E_{R_n(\delta)}(\cdot,\gamma_jn)\Big)=\prod_{j=1}^{\kappa}\mathbb E\big(\mathcal E_{R_n(\delta)}(\cdot,\gamma_jn)\big)+\underset{(\nu,\mu)}{\sum}\sum_{q=1}^{\infty}\sum_D\Big(\frac{e_1}{q}\cdots\frac{e_m}{q}\Big)^nI(D,n,\vecgam,\delta)
\end{align}
where
\begin{align*} 
&I(D,n,\vecgam,\delta)\\
&=V_n^{-2\sum_{j=1}^{\kappa}\gamma_j}\int_{\R^n}\cdots\int_{\R^n}\prod_{j=1}^{\kappa}\bigg(\Big|\sum_{i=1}^m\frac{d_{ij}}{q}\vecx_i\Big|^{-2\gamma_jn}I_{R_n(\delta)}\Big(\sum_{i=1}^m\frac{d_{ij}}{q}\vecx_i\Big)\bigg)\,d\vecx_1\ldots d\vecx_m.
\end{align*}

\begin{prop}\label{propi2}
Let $\kappa\geq2$ and let $\mathcal X_{\kappa}$ denote the set of $\kappa$-admissible matrices $D$, with elements $d_{ij}\in\{0,\pm1\}$, having exactly one nonzero entry in each column. Fix $\delta>0$ and let $\vecgam=(\gamma_1,\ldots,\gamma_{\kappa})$ be such that $\frac{1}{2}<\gamma_1\leq \gamma_2\leq\cdots\leq \gamma_{\kappa}$. Given $D\in\mathcal X_{\kappa}$ let $I_i(D)$, $1\leq i\leq m$, be the set of indices of columns of $D$ that have their nonzero entry in the $i$:th row. Furthermore, set $\Gamma_i(D)=\sum_{j\in I_i(D)}\gamma_j$, $1\leq i\leq m$.  Then
\begin{align*}
 I(D,n,\vecgam,\delta)&=\delta^{m-2\sum_{j=1}^{\kappa}\gamma_j}\prod_{i=1}^{m}\frac{1}{2\Gamma_i(D)-1}
\end{align*}
for all $D\in\mathcal X_{\kappa}$ and all $n\geq1$.
\end{prop}

\begin{Proof}
For $D\in\mathcal X_{\kappa}$ we get
\begin{align*}
I(D,n,\vecgam,\delta)=V_n^{-2\sum_{j=1}^{\kappa}\gamma_j}\prod_{i=1}^{m}\int_{\R^n}|\vecx_i|^{-2\Gamma_i(D)n}I_{R_n(\delta)}(\vecx_i)\,d\vecx_i.
\end{align*}
Using $\sum_{i=1}^{m}\Gamma_i(D)=\sum_{j=1}^{\kappa}\gamma_j$ together with equations \eqref{unitvolume} and \eqref{rnd} yields
\begin{align*}
I(D,n,\vecgam,\delta)&=V_n^{-2\sum_{j=1}^{\kappa}\gamma_j}\prod_{i=1}^{m}\Big(\omega_n\int_{R_n(\delta)}^{\infty}r^{n-2\Gamma_i(D)n-1}\,dr\Big)\\
&=V_n^{-2\sum_{j=1}^{\kappa}\gamma_j}\prod_{i=1}^{m}\Big(\omega_n\frac{R_n(\delta)^{n(1-2\Gamma_i(D))}}{n\big(2\Gamma_i(D)-1\big)}\Big)\\
&=V_n^{m-2\sum_{j=1}^{\kappa}\gamma_j}R_n(\delta)^{n(m-2\sum_{j=1}^{\kappa}\gamma_j)}\prod_{i=1}^{m}\frac{1}{2\Gamma_i(D)-1}\\
&=\delta^{m-2\sum_{j=1}^{\kappa}\gamma_j}\prod_{i=1}^{m}\frac{1}{2\Gamma_i(D)-1}
\end{align*}
for all $n\geq1$.
\end{Proof}

\begin{thm}\label{mixed}
Let $\kappa\geq2$. Fix $\delta>0$ and let $\vecgam=(\gamma_1,\ldots,\gamma_{\kappa})$ be such that $\frac{1}{2}<\gamma_1\leq \gamma_2\leq\cdots\leq \gamma_{\kappa}$. Then
\begin{align*}
\mathbb E\Big(\prod_{j=1}^{\kappa}\mathcal E_{R_n(\delta)}(\cdot,\gamma_jn)\Big)&\to \prod_{j=1}^{\kappa}\frac{\delta^{1-2\gamma_j}}{2\gamma_j-1}+\underset{D\in\mathcal X_{\kappa}}{\sum}\delta^{m-2\sum_{j=1}^{\kappa}\gamma_j}\prod_{i=1}^{m}\frac{1}{2\Gamma_i(D)-1}
\end{align*}
as $n\to\infty$.
\end{thm}

\begin{Proof}
In order to estimate the contribution to \eqref{mix1} from the terms not discussed in Proposition \ref{propi2} we use \eqref{volume} to note that 
\begin{align*}
&I(D,n,\vecgam,\delta)\\
&=\int_{\R^n}\cdots\int_{\R^n}\prod_{j=1}^{\kappa}\bigg(\Big(V_n\Big|\sum_{i=1}^m\frac{d_{ij}}{q}\vecx_i\Big|^{n}\Big)^{-2\gamma_j}I_{R_n(\delta)}\Big(\sum_{i=1}^m\frac{d_{ij}}{q}\vecx_i\Big)\bigg)\,d\vecx_1\ldots d\vecx_m\\
&\ll\int_{\R^n}\cdots\int_{\R^n}\prod_{j=1}^{\kappa}\bigg(\Big(V_n\Big|\sum_{i=1}^m\frac{d_{ij}}{q}\vecx_i\Big|^{n}\Big)^{-2\min_{1\leq j\leq \kappa}\gamma_j}I_{R_n(\delta)}\Big(\sum_{i=1}^m\frac{d_{ij}}{q}\vecx_i\Big)\bigg)\,d\vecx_1\ldots d\vecx_m,
\end{align*}
where the implied constant depends on $\vecgam$ and $\delta$ but not on $n$. Since this last integral is of the form $I(D,n,\min_{1\leq j\leq \kappa}\gamma_j,\delta)$ defined in \eqref{Ek}, it follows from Section \ref{heavy} that the total contribution to \eqref{mix1} from all such terms is exponentially decaying with $n$. 
\end{Proof}

\section{Further results}\label{lastsection}

In this last section we discuss some consequences of the results in Sections \ref{moments} and \ref{joint} as well as the proof of Theorem \ref{curvethm}.

\subsection{The value distribution of $\mathcal E_{R_n(\delta)}(\cdot,\vecc n)$}

First we discuss the random variable $\mathcal E_{R_n(\delta)}(\cdot,cn)$ for $c>\frac{1}{2}$. Let $k\geq1$. We recall from Theorem \ref{higher} that 
\begin{align}\label{highertilde2}
\mathbb E\big(\mathcal E_{R_n(\delta)}(\cdot,cn)^k\big)\to \sum_{D\in\mathcal{D}(k)}2^{k-m}\delta^{m-2kc}\prod_{i=1}^{m}\frac{1}{2k_ic-1}
\end{align}
as $n\to\infty$, where $\mathcal{D}(k)$ is defined as in Lemma \ref{bijection} and $k_i$ is the number of non-zero entries of $D$ in the $i$:th row, $1\leq i\leq m$.

\begin{prop}\label{conv}
Let $k\geq1$, $c>\frac{1}{2}$ and $\delta>0$ be fixed. Then
\begin{align*}
\mathbb E\big(\mathcal E_{R_n(\delta)}(\cdot,cn)^k\big)\to \mathbb E\big(T(c,\delta)^k\big)
\end{align*}
as $n\to\infty$.
\end{prop}

\begin{Proof}
 This is an immediate consequence of \eqref{Tmoment}, \eqref{highertilde2} and Lemma \ref{bijection}.
\end{Proof}

\begin{cor}
Let $c>\frac{1}{2}$ and $\delta>0$ be fixed. Then $\mathcal E_{R_n(\delta)}(\cdot,cn)$ converges in distribution to $T(c,\delta)$ as $n\to\infty$.
\end{cor}

\begin{Proof}
Since the distribution of $T(c,\delta)$ is uniquely determined by its moments by Lemma \ref{determ} and
\begin{align*}
 \lim_{n\to\infty}\mathbb E\big(\mathcal E_{R_n(\delta)}(\cdot,cn)^k\big)=\mathbb E\big(T(c,\delta)^k\big)
\end{align*}
for all $k\geq1$ by Proposition \ref{conv}, the corollary follows from \cite[Thm.\ 30.2]{billing}.
\end{Proof}

We now continue with the general situation. Let $m\geq1$ and fix $\vecc=(c_1,\ldots,c_m)$ satisfying $\frac{1}{2}<c_1< c_2<\cdots<c_m$. We will discuss the random vector $\mathcal E_{R_n(\delta)}(\cdot,\vecc n)$. For $\kappa\geq2$  and fixed $\gamma_j\in\{c_1,\ldots,c_m\}$, $1\leq j\leq\kappa$, satisfying $\gamma_1\leq \gamma_2\leq\cdots\leq \gamma_{\kappa}$, we recall from Section \ref{joint} that 
\begin{align}\label{mixedtilde2}
\mathbb E\Big(\prod_{j=1}^{\kappa}\mathcal E_{R_n(\delta)}(\cdot,\gamma_jn)\Big)&\to\sum_{D\in\mathcal{D}(\kappa)}2^{\kappa-m}\delta^{m-2\sum_{j=1}^{\kappa}\gamma_j}\prod_{i=1}^{m}\frac{1}{2\Gamma_i(D)-1}
\end{align}
as $n\to\infty$, where $\Gamma_i(D)$, $1\leq i\leq m$, is defined as in Proposition \ref{propi2}.

\begin{prop}\label{mixedconv}
Fix $\delta>0$. Let $\kappa\geq2$ and $\vecgam=(\gamma_1,\ldots,\gamma_{\kappa})$ be such that $\frac{1}{2}<\gamma_1\leq \gamma_2\leq\cdots\leq \gamma_{\kappa}$. Then
\begin{align*}
\mathbb E\Big(\prod_{j=1}^{\kappa}\mathcal E_{R_n(\delta)}(\cdot,\gamma_jn)\Big)\to \mathbb E\Big(\prod_{j=1}^{\kappa}T(\gamma_j,\delta)\Big)
\end{align*}
as $n\to\infty$.
\end{prop}

\begin{Proof}
This is an immediate consequence of \eqref{Tmixedmoment}, \eqref{mixedtilde2} and Lemma \ref{bijection}.
\end{Proof}

\begin{remark}
Propositions \ref{conv} and \ref{mixedconv} together complete the proof of Theorem \ref{momentthm}. 
\end{remark}

\begin{cor}\label{mixedconvergence}
Fix $\delta>0$. Let $m\geq1$ and let $\vecc=(c_1,\ldots,c_{m})$ be such that $\frac{1}{2}<c_1< c_2<\cdots< c_{m}$.
Then the random vector $\mathcal E_{R_n(\delta)}(\cdot,\vecc n)$ converges in distribution to $T(\vecc,\delta)$ as $n\to\infty$.
\end{cor}

\begin{Proof}
Since the distribution of $T(\vecc,\delta)$ is uniquely determined by its joint moments by Lemma \ref{multideterm} and
the joint moments of $\mathcal E_{R_n(\delta)}(\cdot,\vecc n)$ converge to those of $T(\vecc,\delta)$ as $n\to\infty$
by Proposition \ref{conv} and Proposition \ref{mixedconv}, the corollary follows from a (in principle word by word) generalization of \cite[Thm.\ 30.2]{billing}.
\end{Proof}

\subsection{An alternative proof of Theorem \ref{mainthm}}\label{altsec}

Let $\vecc=(c_1,\ldots,c_m)$ and recall the definition of the random vectors $\mathcal E_n(\cdot,\vecc n)$, $\mathcal E_{R_n(\delta)}(\cdot,\vecc n)$, $T(\vecc)$ and $T(\vecc,\delta)$ from \eqref{Evec}, \eqref{Edvec}, \eqref{Tvec} and \eqref{Tdvec}. We call their distributions $\mu_{\mathcal E_n(\cdot,\vecc n)}$, $\mu_{\mathcal E_{R_n(\delta)}(\cdot,\vecc n)}$, $\mu_{T(\vecc)}$ and $\mu_{T(\vecc,\delta)}$ respectively. The idea of the proof is to let $\delta\to0$ in Corollary \ref{mixedconvergence}. The details are as follows.

The result in Corollary \ref{mixedconvergence} can be restated, using the Lévy-Prohorov metric (cf.\ Section \ref{mainsec}), as for every $\ve>0$ there exists $N_{\delta,\ve}\in\Z^+$ such that for every $n\geq N_{\delta,\ve}$ we have
\begin{align}\label{firstalt}
\pi(\mu_{\mathcal E_{R_n(\delta)}(\cdot,\vecc n)},\mu_{T(\vecc,\delta)})<\ve. 
\end{align}
We also note that, for fixed $\delta>0$, it follows from \eqref{secondalt} that for every Borel set $A\subseteq\R^m$ we have
\begin{align*}
\mu_{\mathcal E_{R_n(\delta)}(\cdot,\vecc n)}(A)\leq\mu_{\mathcal E_n(\cdot,\vecc n)}(A)+\frac{\delta}{2}, 
\end{align*}
which in turn implies that
\begin{align}\label{thirdalt}
\pi(\mu_{\mathcal E_{R_n(\delta)}(\cdot,\vecc n)},\mu_{\mathcal E_n(\cdot,\vecc n)})\leq\frac{\delta}{2}. 
\end{align}
We stress that inequality \eqref{thirdalt} holds independently of $n\geq1$.

In order to obtain a similar estimate on the Lévy-Prohorov distance between $\mu_{T(\vecc)}$ and $\mu_{T(\vecc,\delta)}$ we observe that
\begin{align*}
\text{Prob}\big(T(\vecc,\delta)\neq T(\vecc)\big)=\text{Prob}\big(T_1\leq \delta\big)=\text{Prob}\big(N(\delta)\geq1\big)=1-e^{-\delta/2}<\frac{\delta}{2}. 
\end{align*}
Here we have used the fact that $N(\delta)$ is Poisson distributed with expectation value $\frac{1}{2}\delta$. Hence, for every Borel set $A\subseteq\R^m$, we have  
\begin{align*}
\mu_{T(\vecc,\delta)}(A)<\mu_{T(\vecc)}(A)+\frac{\delta}{2}, 
\end{align*}
which implies that
\begin{align}\label{fourthalt}
\pi(\mu_{T(\vecc,\delta)},\mu_{T(\vecc)})\leq\frac{\delta}{2}. 
\end{align}

To finish the proof let $\ve>0$ be given and let $0<\delta\leq\ve$. Now \eqref{firstalt}, \eqref{thirdalt} and \eqref{fourthalt} yield
\begin{align*}
\pi(\mu_{\mathcal E_n(\cdot,\vecc n)},\mu_{T(\vecc)})&\leq\pi(\mu_{\mathcal E_n(\cdot,\vecc n)},\mu_{\mathcal E_{R_n(\delta)}(\cdot,\vecc n)})\\
&+
\pi(\mu_{\mathcal E_{R_n(\delta)}(\cdot,\vecc n)},\mu_{T(\vecc,\delta)})+
\pi(\mu_{T(\vecc,\delta)},\mu_{T(\vecc)})<2\ve 
\end{align*}
for all $n\geq N_{\delta,\ve}$. We conclude that $\mu_{\mathcal E_n(\cdot,\vecc n)}$ converges in the metric $\pi$ to $\mu_{T(\vecc)}$ as $n\to\infty$. This means that $\mathcal E_n(\cdot,\vecc n)$ converges in distribution to $T(\vecc)$ as $n\to\infty$.

\subsection{Random functions and the proof of Theorem \ref{curvethm}} 

Fix $\frac{1}{2}<A<B$. In this section we study
\begin{align*}
\widehat{\mathcal E_n}(\cdot,c):[A,B]\to \R,\,\, c\mapsto \mathcal E_n(\cdot,cn)
\end{align*}
as a random function in $C[A,B]$. Our aim is to understand the limit distribution of this random function as $n\to\infty$. We will derive Theorem \ref{curvethm} as a consequence of Proposition \ref{convex} below. First we prove the following elementary lemma.

\begin{lem}\label{easyconvex}
Let $0<\delta\leq1$ and $\ve>0$. Let $f\in C[0,1]$ be convex and suppose that there exist $0\leq t<s\leq\min(1,t+\delta)$ such that $|f(s)-f(t)|\geq\ve$. Then either $|f(\sfrac{\delta}{2})-f(0)|\geq\frac{\ve}{2}$ or $|f(1)-f(1-\sfrac{\delta}{2})|\geq\frac{\ve}{2}$.
\end{lem}

\begin{Proof}
Since $f$ is convex we can without loss of generality assume that $t=0$ or $s=1$. Say that $t=0$ (the case $s=1$ is treated analogously). If $f$ is monotone on $[0,\delta]$ the result is obvious.
If $f$ is not monotone on $[0,\delta]$ and $|f(\sfrac{\delta}{2})-f(0)|<\frac{\ve}{2}$, then  
\begin{align*}
\frac{\ve}{2}\leq f(\delta)-f(\sfrac{\delta}{2})\leq f(1)-f(1-\sfrac{\delta}{2})
\end{align*}
by the convexity of $f$.
\end{Proof}

\begin{prop}\label{convex}
Let $P_n$ and $P$ be Borel probability measures on $C[0,1]$. Assume that for every $n$, $P_n$-almost all $f\in C[0,1]$ are convex. If all the finite-dimensional distributions of $P_n$ converge weakly to those of $P$, then $P_n$ converge weakly to $P$. 
\end{prop}

\begin{Proof}
By \cite[Thm.\ 7.1]{billconv} it suffices to prove that the sequence $\{P_n\}$ is tight. First of all, since $\lim_{a\to\infty}P\{f:|f(0)|\geq a\}=0$ (cf.\ e.g.\ \cite[Thm.\ 1.19(e)]{rudin}) and $\limsup_{n\to\infty}P_n\{f:|f(0)|\geq a\}\leq P\{f:|f(0)|\geq a\}$ for every $a$ (cf.\ \cite[Thm.\ 2.1(iii)]{billconv}), we note that for each $\eta>0$ there exist $a>0$ and $n_0\in\Z^+$ such that 
\begin{align}\label{COND1}
P_n\{f:|f(0)|\geq a\}\leq\eta \qquad\text{for all }\: n\geq n_0.
\end{align}

Next let $\ve,\eta>0$ be given. Choose $k\in\Z^+$ such that
\begin{align*}
P\Big\{f\in C[0,1]:\big|f(\sfrac{1}{2k})-f(0)\big|\geq\frac{\ve}{2} \text{ or } \big|f(1)-f(\sfrac{2k-1}{2k})\big|\geq\frac{\ve}{2}\Big\}<\eta 
\end{align*}
(as is possible by basic measure theory, cf.\ \cite[Thm.\ 1.19(e)]{rudin}). Then, by \cite[Thm.\ 2.1(iii)]{billconv}, there exists $n_0\in\Z^+$ such that 
\begin{align}\label{<eta}
P_n\Big\{f\in C[0,1]:\big|f(\sfrac{1}{2k})-f(0)\big|\geq\frac{\ve}{2} \text{ or } \big|f(1)-f(\sfrac{2k-1}{2k})\big|\geq\frac{\ve}{2}\Big\}\leq\eta 
\end{align}
for all $n\geq n_0$. For every $n\geq n_0$, $P_n$-almost every $f\in C[0,1]$ which does not lie in the set in \eqref{<eta} is convex. For every such convex $f$ it follows from Lemma \ref{easyconvex} that $|f(s)-f(t)|<\ve$ for all $0\leq t\leq s\leq \min(1,t+k^{-1})$. We conclude that for every $n\geq n_0$ we have 
\begin{align}\label{COND2}
P_n\Big\{f\in C[0,1]:\sup_{|s-t|\leq k^{-1}}|f(s)-f(t)|>\ve\Big\}\leq\eta. 
\end{align}
In view of \eqref{COND1} and \eqref{COND2} the sequence $\{P_n\}$ is tight (cf.\ \cite[Thm.\ 7.3]{billconv}), and this completes the proof of the proposition.
\end{Proof}

\begin{Proof2}
Since, for each fixed $L\in X_n$, $\widehat{\mathcal E_n}(L,c)$ is convex, Theorem \ref{curvethm} follows from Proposition \ref{convex} together with Theorem \ref{mainthm}.
\end{Proof2}

Finally we extend Theorem \ref{curvethm} to semi-infinite intervals.

\begin{cor}
For each $n\in\Z^+$ and any fixed $A>\frac{1}{2}$ consider
\begin{align*}
 c\mapsto V_n^{-2c}E_n(\cdot,cn)
\end{align*}
as a random function in $C[A,\infty)$. The distribution of this random function converges to the distribution of
\begin{align*}
 c\mapsto 2\sum_{j=1}^{\infty}T_j^{-2c}
\end{align*}
as $n\to\infty$.
\end{cor}

\begin{Proof}
This follows from \cite[Thm.\ 5]{whitt} and Theorem \ref{curvethm}. 
\end{Proof}

\subsubsection*{Acknowledgement} The author is most grateful to Andreas Strömbergsson for many helpful discussions on this work.

\end{document}